%% file: main.tex
\documentclass[preprint,12pt]{elsarticle}
\usepackage{amsmath,amssymb,amsthm,mathtools,bm,esint,mathrsfs}
\usepackage{booktabs,array}
\usepackage{graphicx}
\usepackage{tikz}
\usetikzlibrary{arrows.meta,patterns}
\usepackage{enumitem}
\usepackage{microtype}
\usepackage{needspace}
\usepackage[above,below]{placeins}
\usepackage{float}
\usepackage{comment}
\usepackage{xcolor}
\usepackage[hidelinks]{hyperref}
\usepackage[nameinlink,noabbrev]{cleveref}
\journal{Journal of the Mechanics and Physics of Solids}

\newtheorem{theorem}{Theorem}[section]
\newtheorem{proposition}[theorem]{Proposition}

\theoremstyle{definition}

\theoremstyle{remark}
\newtheorem{remark}[theorem]{Remark}

\newcommand{\dd}{\,\mathrm d}
\newcommand{\E}{\mathcal E}
\newcommand{\R}{\mathbb{R}}

\newcommand{\Fu}{\mathcal G_\varepsilon}
\newcommand{\Et}{\widetilde{\mathcal E}}
\newcommand{\sym}{\mathrm{sym}}
\newcommand{\He}{\mathcal{H}}
\newcommand{\Aeps}{\mathcal{A}_{\varepsilon}}

\newcommand{\Fsplit}{\mathcal F_\varepsilon^{\mathrm{split}}}
\newcommand{\Fsplith}{\mathcal F_{\varepsilon,h}^{\mathrm{split}}}
\newcommand{\Fustar}{\mathcal G_{\varepsilon,*}}
\newcommand{\Aepsstar}{\mathcal{A}_{\varepsilon}^{*}}
\newcommand{\sech}{\operatorname{sech}}
\newcommand{\norm}[1]{\left\lVert#1\right\rVert}
\newcommand{\pair}[2]{\left\langle #1,#2\right\rangle}
\newcommand{\revision}[1]{\textcolor{black}{#1}}

\begin{document}

\begin{frontmatter}

\title{Nonlocal regularization of sharp-crack energetics:
a variational bridge from Griffith geometry to emergent constitutive degradation}

\author[mathematics]{Qing Cheng}
\ead{qingcheng@tongji.edu.cn}
\author[mathematics]{Yuqi Sun}
\ead{yuqisun@tongji.edu.cn}
\author[mathematics]{Xuejun Xu}
\ead{xuxj@tongji.edu.cn}

\address[mathematics]{School of Mathematical Sciences, Key Laboratory of
Intelligent Computing and Applications (Ministry of Education),
Tongji University, Shanghai, 200092, China}


\begin{abstract}
We develop a variational diffuse-crack model based on nonlocal
averaging of tensile strain energy. A bounded potential acts on the scaled
average, while the compressive energy remains local. \revision{The initial slope preserves the averaged tensile strain-energy density
to leading order.}
The first variation yields a constitutive stress in which the tensile
contribution is multiplied by the adjoint average of the derivative of the
saturation potential. The resulting degradation is therefore determined by
the stored-energy functional, without introducing an independent fracture
field or prescribing a separate degradation law. For fixed smooth displacement fields with bounded strain, the energy
density converges to the classical elastic strain-energy density as the
averaging radius tends to zero. The Griffith surface-energy coefficient
is recovered for planar opening transitions.
A structure-preserving scalar auxiliary variable (SAV) scheme satisfies an
exact discrete balance between the \revision{modified-energy increment} and external work.
\revision{We compare SAV and velocity-Verlet solutions for dynamic crack
branching and Kalthoff--Winkler impact.}
\end{abstract}

\begin{keyword}
nonlocal regularization \sep variational fracture
\sep tension--compression splitting
\sep energy saturation \sep emergent constitutive degradation\sep structure-preserving time integration
\end{keyword}

\end{frontmatter}

\section{Introduction}
\label{sec:introduction}
\noindent Fracture couples solid deformation to the creation of an internal surface.
Griffith's energetic description balances elastic-energy release against the energy required to create new crack surfaces
\citep{Griffith1921}. Variational fracture treats
displacement and crack geometry as unknowns
\citep{FrancfortMarigo1998,BourdinFrancfortMarigo2008}, allowing crack paths
to emerge from energetic competition. Computing this evolution requires a numerical representation of
displacement jumps across the growing crack surface. Enriched finite
elements accommodate these jumps without a crack-conforming mesh
\citep{MoesDolbowBelytschko1999}, but still require a separate representation
of the crack geometry that is updated during propagation.

Diffusive fracture methods instead approximate the crack-surface
energy by a volume integral of a fracture-energy density concentrated in a
narrow region around the crack. Building on elliptic approximations of free-discontinuity
functionals \citep{AmbrosioTortorelli1990}, phase-field fracture couples
displacement to a scalar crack field \citep{BourdinFrancfortMarigo2000}.
In this framework, a crack-surface density defines the regularized
fracture energy, while a degradation function couples the elastic energy
to the crack field. Cohesive phase-field extensions use these functions to
represent traction--separation laws
\citep{Wu2017,ContiFocardiIurlano2016}.
Phase-field models have been applied to a broad range of fracture
problems, including dynamic fracture
\revision{\citep{BordenEtAl2012,ChengSunXu2026}}, ductile fracture
in elastoplastic solids \citep{AmbatiGerasimovDeLorenzis2015Ductile},
finite-strain fracture in hyperelastic solids \citep{MieheSchanzel2014},
and multiphysics fracture
\citep{MieheSchanzelUlmer2015,MieheEtAl2016Chemo,MartinezPanedaGolahmarNiordson2018}.
\revision{Across these applications, phase-field models avoid explicit
crack-surface tracking by solving for an additional fracture field coupled
to displacement.}

Nonlocal and gradient-enhanced damage models introduce spatial scales to
regularize softening and localization
\citep{PijaudierCabotBazant1987,PeerlingsEtAl1996}, and variational gradient
damage connects these processes to brittle-fracture energetics
\citep{PhamAmorMarigoMaurini2011}. \revision{Other routes include virtual and geometric internal-bond models
\citep{GaoKlein1998,SunYu2025}, peridynamics
\citep{Silling2000,SillingEtAl2007}, and nonlocal geometric fracture theory
\citep{YuSun2023}.}
Eliminating a local damage variable can
also produce a saturating effective energy \citep{RenEtAl2025}. \revision{The elimination of a local damage variable suggests a formulation
with displacement as its sole primary variable. The aim here is to combine
this feature with nonlocal spatial regularization in a variational
diffuse-crack model grounded in sharp-crack energetics.}

Nonlocal approximations of free-discontinuity energies connect spatial
averages of gradient- or strain-dependent energy densities to limiting
functionals with surface-energy contributions on the jump set. \citet{BraidesDalMaso1997} established a nonlocal
$\Gamma$-approximation of the scalar Mumford--Shah functional.
For unsplit elastic energies based on linearized strain, Griffith-type
$\Gamma$-limits have been established for nonlocal energies in both
continuum and finite-element settings
\citep{LussardiNegri2007,ScillaSolombrino2021}.

\revision{We develop a variational diffuse-crack model for dynamic fracture based
on nonlocal averaging of tensile strain energy.} A bounded potential acts on the scaled tensile-energy average,
while the compressive energy remains local. Tensile stress degradation \revision{emerges} from the exact first variation, with displacement as the sole
primary field. We further construct a structure-preserving time
discretization for the resulting dynamics. The specific contributions
are as follows:
\begin{enumerate}[label=(\roman*),leftmargin=*,itemsep=0.4em]
 \item We specify the potential $\Phi$ to recover the classical
 elastic strain-energy density for fixed smooth displacement fields as
 the averaging radius tends to zero, and the Griffith fracture energy per
 unit area for planar opening. With its initial slope and limiting value
 fixed, the shape of $\Phi$ determines the homogeneous tensile strength
 at a fixed averaging radius.
 \item We introduce tension--compression splitting into the nonlocal
 energy, saturating only the averaged tensile contribution and retaining
 local compressive energy. Constitutive degradation emerges
 from the exact first variation of the energy, without an independent
 fracture field or a separately prescribed degradation law.

 \item We construct a linearly implicit, structure-preserving
 scalar auxiliary variable (SAV) scheme
 that satisfies an exact discrete balance between the modified-energy
increment and the external work.
\end{enumerate}

\revision{We compare the SAV and velocity-Verlet solutions for two standard planar
dynamic-fracture benchmarks: dynamic crack branching and
Kalthoff--Winkler impact.}
The split--unsplit comparison indicates that tension--compression
splitting provides a more physically plausible description of mixed-mode
dynamic cracking.

Section~\ref{sec:variational-governing-sharp-interface}
introduces the model with tension--compression splitting and derives
its constitutive stress and governing equations.
Section~\ref{sec:mechanical-meaning} examines the elastic and
opening-crack limits, emergent constitutive degradation and homogeneous
strength.
Section~\ref{sec:sav} develops the SAV scheme and its
discrete energy law. Section~\ref{sec:numerics} presents the spatial
discretization and the two dynamic benchmarks.

\section{Variational formulation and governing equations}
\label{sec:variational-governing-sharp-interface}
\noindent \revision{We regularize the Griffith energy by applying a bounded
potential to the scaled nonlocal average of tensile strain energy and retaining
the compressive energy locally. We then derive the constitutive stress
and governing equations.}

\subsection{The Griffith energy of sharp cracks}
\noindent Let $\Omega\subset\R^d$ be a bounded open set, $d\geq2$, and let
$\mathbf{u}:\Omega\to\R^d$ denote the displacement field. The infinitesimal strain tensor is defined by
\begin{equation}
  \mathbf{e}(\mathbf{u})=\sym\nabla \mathbf{u}
  =\frac{1}{2}\bigl(\nabla\mathbf{u}+(\nabla\mathbf{u})^{\mathsf T}\bigr).
\end{equation}
The jump set $J_{\mathbf u}$ of the displacement field is defined through
its one-sided approximate limits.  For $\mathbf{x}\in\Omega$, $\rho>0$, and
$\bm{\nu}\in\mathbb S^{d-1}$, let
\begin{equation}
  B_\rho^\pm(\mathbf{x},\bm{\nu})
  :=\bigl\{\mathbf{y}\in B_\rho(\mathbf{x}):\ \pm(\mathbf{y}-\mathbf{x})\cdot\bm{\nu}>0\bigr\}.
  \label{eq:half-balls}
\end{equation}
Then
\begin{equation}
\begin{aligned}
  J_{\mathbf u}:=\Bigl\{&\mathbf{x}\in\Omega:\;
  \exists\,\mathbf u^\pm(\mathbf{x})\in\R^d,
  \ \bm{\nu}_{\mathbf u}(\mathbf{x})\in\mathbb S^{d-1},\, \mathbf u^+(\mathbf{x})\neq\mathbf u^-(\mathbf{x}),\\[-0.2em]
  &\lim_{\rho\downarrow0}
  \fint_{B_\rho^\pm(\mathbf{x},\bm{\nu}_{\mathbf u}(\mathbf{x}))}
  \bigl(1\wedge\lvert\mathbf u(\mathbf{y})-\mathbf u^\pm(\mathbf{x})\rvert\bigr)\dd\mathbf{y}=0\Bigr\}.
  \label{eq:jump-set}
\end{aligned}
\end{equation}
where $a\wedge b:=\min\{a,b\}$ denotes the minimum. The total potential energy, including the Griffith fracture energy, is defined by
\begin{equation}
  \Pi_0[\mathbf{u}]
  =\int_{\Omega\setminus J_{\mathbf{u}}} W(\mathbf{e}(\mathbf{u}))\dd\mathbf{x}
   +G_c\He^{d-1}(J_{\mathbf u})
   -\ell_{\mathrm{ext}}[\mathbf{u}], 
  \label{eq:griffith}
\end{equation}
where $W$ is the elastic strain-energy density, $G_c>0$ is the critical
energy-release rate, $\ell_{\mathrm{ext}}[\mathbf{u}]$ is the external work functional, and $\He^{d-1}(J_{\mathbf u})$ is the Hausdorff measure of the sharp-crack interface.

The natural setting for \eqref{eq:griffith} is the space of generalized
special functions of bounded deformation, $GSBD^p(\Omega)$, which admits
displacement jumps \citep{DalMaso2013}. The main difficulty is that the crack
geometry $J_{\mathbf u}$ is unknown
\citep{Griffith1921,FrancfortMarigo1998}.

\subsection{Nonlocal regularization with a tension--compression split}
\noindent To retain compressive stiffness in the regularized model,
we first separate the tensile and compressive contributions to the elastic
strain energy. For an isotropic linearly elastic solid, the strain-energy
density is
\begin{equation}
  W(\mathbf e)=\frac{\lambda}{2}(\operatorname{tr}\mathbf e)^2
                 +\mu\,\mathbf e:\mathbf e,
\end{equation}
where $\lambda\geq0$ and $\mu>0$ are Lamé parameters. Following the strain-spectral decomposition
\citep{MieheHofackerWelschinger2010}, define
\begin{equation}
 \mathbf e=\sum_{\alpha=1}^{d}\eta_\alpha
       \mathbf n_\alpha\otimes\mathbf n_\alpha,\quad
  \mathbf e^\pm=\sum_{\alpha=1}^{d}\langle\eta_\alpha\rangle_\pm
       \mathbf n_\alpha\otimes\mathbf n_\alpha.
\label{eq:spectral-energy-split}
\end{equation}
where $\eta_\alpha$ and $\mathbf n_\alpha$ are the eigenvalues and
corresponding orthonormal eigenvectors of $\mathbf e$, respectively, and
$\langle r\rangle_+=\max(r,0)$ and $\langle r\rangle_-=\min(r,0)$.  The elastic strain-energy density then admits the decomposition
\begin{equation}
W(\mathbf e)=W^+(\mathbf e)+W^-(\mathbf e),
\end{equation}
where $  W^\pm(\mathbf e)=\frac{\lambda}{2}
\langle\operatorname{tr}\mathbf e\rangle_\pm^2
+\mu\,\mathbf e^\pm:\mathbf e^\pm$. 

For a scalar field $f$, define the nonlocal averaging operator by
\begin{equation}
  (\Aeps f)(\mathbf x)
  =\frac{1}{m_\varepsilon(\mathbf x)}
    \int_{B_\varepsilon(\mathbf x)\cap\Omega}f(\mathbf y)\dd\mathbf y,
  \label{eq:average}
\end{equation}
where $B_\varepsilon(\mathbf x)$ is the open ball centered at $\mathbf x$ with radius $\varepsilon$ and $m_\varepsilon(\mathbf x)=|B_\varepsilon(\mathbf x)\cap\Omega|$. \revision{To regularize tensile fracture, we apply a bounded saturation potential
to the scaled nonlocal average of tensile strain energy.}
The lower-semicontinuous, nondecreasing potential
$\Phi:[0,\infty)\to[0,\infty)$ satisfies
\begin{equation}
  \Phi(0)=0,\quad
  \lim_{s\downarrow0}\frac{\Phi(s)}{s}=1,\quad
  \lim_{s\to\infty}\Phi(s)=\frac{G_c}{2}.
  \label{eq:saturation-normalization}
\end{equation}
The small-$s$ condition implies $\Phi(s)=s+o(s)$ and preserves the
leading-order tensile elastic contribution. The saturation value $G_c/2$
\revision{recovers the Griffith surface energy in the planar-opening limit.
We use the exponential saturation potential}
\begin{equation}
	\Phi_{\exp}(s)
	=\frac{G_c}{2}\left(1-\exp\left[-\frac{2s}{G_c}\right]\right).
	\label{eq:exponential-potential}
\end{equation}
\revision{\ref{app:saturation-laws} lists other saturation potentials.}

\begin{remark}
\revision{In three dimensions,}
$W$ has units J/m$^3$, so $s=\varepsilon\Aeps W$ has units J/m$^2$,
the same as $G_c$. Thus $\Phi$ also has units J/m$^2$ and
$\Phi(s)/\varepsilon$ is an energy density.
\end{remark}

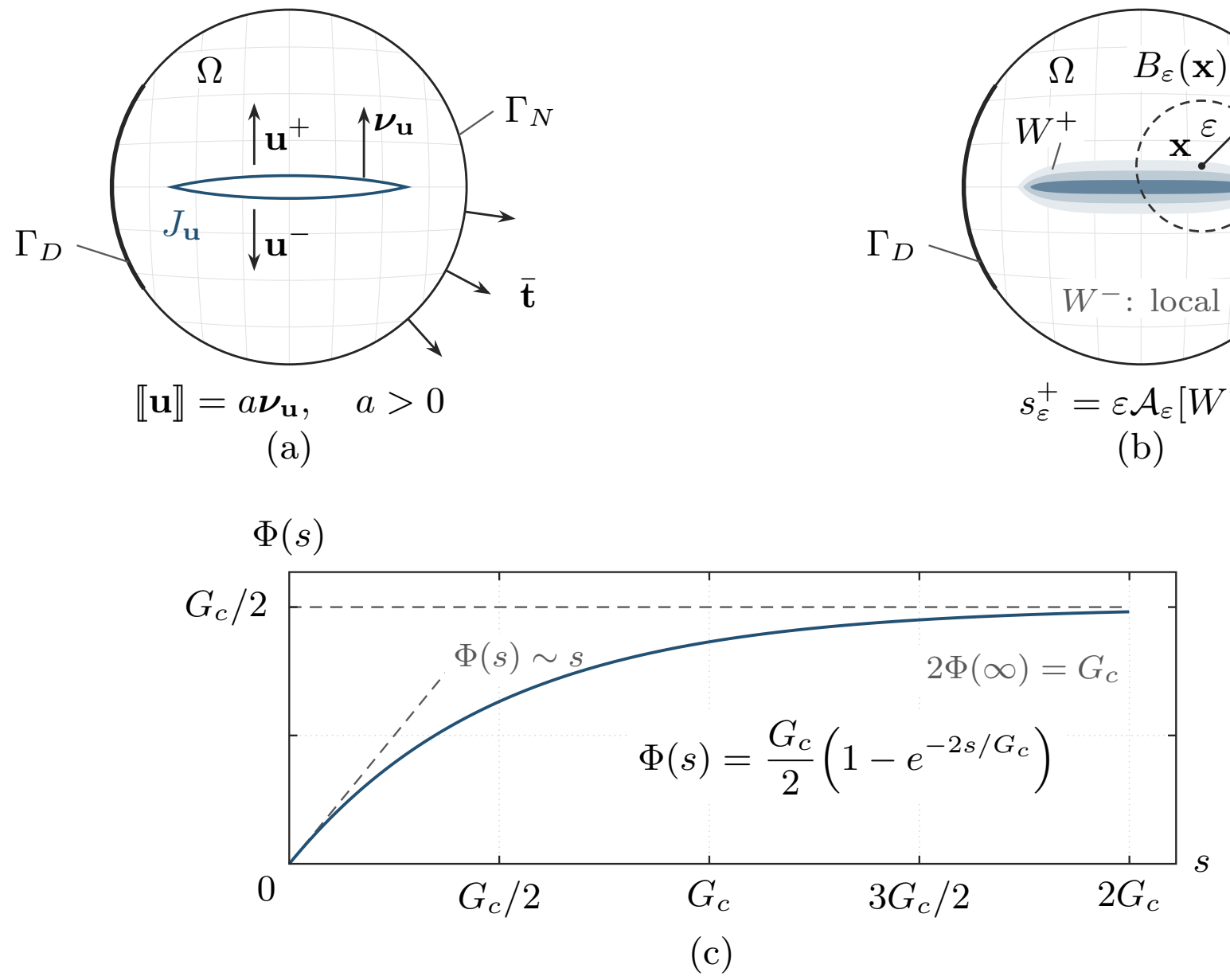
\begin{figure}[!htbp]
\centering
\resizebox{\linewidth}{!}{\input{figures/split_functional_schematic.tex}}
\caption{Nonlocal regularization of fracture in $\Omega$, with
prescribed displacements on $\Gamma_D$ and tractions on $\Gamma_N$:
(a) a sharp opening crack;
(b) diffuse tensile localization and averaging over $B_\varepsilon(\mathbf x)$,
with $W^-$ retained locally; (c) the exponential saturation potential,
with unit initial slope and plateau $G_c/2$.
The opening and localization widths are schematic; $\varepsilon$ denotes
the averaging radius.}
\label{fig:split-functional-schematic}
\end{figure}

\revision{The split stored-energy functional is}
\begin{equation}
\begin{aligned}
  \mathcal F_\varepsilon^{\mathrm{split}}[\mathbf u]
     &=\int_\Omega W^-(\mathbf e(\mathbf u))\dd\mathbf x
       +\frac1\varepsilon\int_\Omega
          \Phi(s_\varepsilon^+[\mathbf u])\dd\mathbf x.
\end{aligned}
  \label{eq:split-functional}
\end{equation}
Here, $s_\varepsilon^+[\mathbf u]=\varepsilon\Aeps[W^+(\mathbf e(\mathbf u))]$ is the scaled nonlocal
average of the tensile strain-energy density. Figure~\ref{fig:split-functional-schematic} illustrates the nonlocal
regularization of an opening crack and the associated saturation potential.

Let $\Omega$ be a fixed Lipschitz reference domain, with disjoint Dirichlet
and Neumann boundary parts $\Gamma_D$ and $\Gamma_N$. For $t\in[0,T]$, impose $\mathbf u=\bar{\mathbf u}(t)$ on $\Gamma_D$ and define
\begin{equation}
  \Pi_\varepsilon^{\mathrm{split}}[\mathbf u;t]
    =\mathcal F_\varepsilon^{\mathrm{split}}[\mathbf u]
      -\ell_{\mathrm{ext}}(t)[\mathbf u],
\end{equation}
where the external work functional associated with the prescribed body
force per unit reference volume $\mathbf b$ and traction $\bar{\mathbf t}$ is
\begin{equation}
  \ell_{\mathrm{ext}}(t)[\mathbf{u}]
  =\int_\Omega \mathbf{b}(\mathbf x,t)\cdot \mathbf{u}\dd\mathbf{x}
   +\int_{\Gamma_N}\bar{\mathbf{t}}(\mathbf x,t)\cdot\mathbf{u}\dd a.
  \label{eq:dead-load-potential}
\end{equation}

\subsection{Energy variations}
\noindent The constitutive stress follows from the first variation
of the stored-energy functional \eqref{eq:split-functional}. To identify
the stress conjugate to the strain variation, we first introduce the
adjoint of the nonlocal averaging operator.
The integral kernel of \eqref{eq:average} is
\begin{equation}
K_\varepsilon(\mathbf{x},\mathbf{y})
=\frac{\bm{1}_{\{|\mathbf{x}-\mathbf{y}|<\varepsilon\}}}
{m_\varepsilon(\mathbf{x})}.
\end{equation}
Accordingly, the adjoint of the averaging operator is
\begin{equation}
  (\Aepsstar g)(\mathbf{y})
  =\int_{B_\varepsilon(\mathbf{y})\cap\Omega}
    \frac{g(\mathbf{x})}{m_\varepsilon(\mathbf{x})}\dd\mathbf{x}.
  \label{eq:adjoint-average}
\end{equation}
In general, \(\Aepsstar\neq\Aeps\) because the normalization depends on the
centre \(\mathbf{x}\) of the truncated ball.

\begin{proposition}[First variation and constitutive stress]
\label{thm:first-variation}
Let \(\mathbf{u}\) be an admissible displacement and \(\mathbf{v}\) an admissible variation,
with \(\mathbf{u},\mathbf{v}\in W^{1,p}(\Omega;\R^d)\), $p>1$, and let
\(\mathbf{v}=\mathbf{0}\) on \(\Gamma_D\).
Assume that
\(W^\pm:\R_{\mathrm{sym}}^{d\times d}\to[0,\infty)\) are $C^1$ and
\begin{equation}
  \lvert D_{\mathbf{E}}W^\pm(\mathbf{E})\rvert
  \leq C_W(1+\lvert \mathbf{E}\rvert^{p-1})
  \quad\text{for all }\mathbf{E}\in\R_{\mathrm{sym}}^{d\times d},
  \label{eq:stress-growth}
\end{equation}
that \(\Phi\in C^1([0,\infty))\) has bounded derivative, and that
\(\ell_{\mathrm{ext}}\) is differentiable at \(\mathbf{u}\).  Define
\begin{equation}
  \boldsymbol{\sigma}^\pm(\mathbf{u})
  =D_{\mathbf{E}}W^\pm(\mathbf{e}(\mathbf{u})),
  \quad
  q_\varepsilon^+[\mathbf{u}](\mathbf{x})
  =\Phi'\bigl(s_\varepsilon^+[\mathbf{u}](\mathbf{x})\bigr),
  \label{eq:q-def}
\end{equation}
and
\begin{equation}
  a_\varepsilon^+[\mathbf{u}](\mathbf{y})
  =(\Aepsstar q_\varepsilon^+[\mathbf{u}])(\mathbf{y}).
  \label{eq:a-def}
\end{equation}
Then
\begin{equation}
  D\Pi_\varepsilon^{\mathrm{split}}[\mathbf{u};t](\mathbf{v})
  =\int_\Omega \bigl(\boldsymbol{\sigma}^-(\mathbf u)
      +a_\varepsilon^+[\mathbf u]\boldsymbol{\sigma}^+(\mathbf u)\bigr)
      :\mathbf{e}(\mathbf{v})\dd\mathbf{y}
   -\ell_{\mathrm{ext}}(t)[\mathbf{v}],
  \label{eq:first-variation}
\end{equation}
and the variationally consistent stress is
\begin{equation}
  \boldsymbol{\sigma}_\varepsilon^{\mathrm{split}}
    =\boldsymbol{\sigma}^-+a_\varepsilon^+\boldsymbol{\sigma}^+,
  \quad
  \boldsymbol{\sigma}^\pm
    =\lambda\langle\operatorname{tr}\mathbf e\rangle_\pm\mathbf I
       +2\mu\mathbf e^\pm
  \label{eq:split-stress}
\end{equation}
for \eqref{eq:spectral-energy-split}.
\end{proposition}

\begin{proof}
Differentiating $s_\varepsilon^+[\mathbf u]$ in the direction $\mathbf v$ gives
\begin{equation}
  Ds_\varepsilon^+[\mathbf{u}](\mathbf{v})
  =\varepsilon\Aeps\bigl[
    \boldsymbol{\sigma}^+(\mathbf{u}):\mathbf{e}(\mathbf{v})\bigr].
\end{equation}
Applying the chain rule and then the adjoint relation yields
\begin{align}
  D\mathcal F_\varepsilon^{\mathrm{split}}[\mathbf{u}](\mathbf{v})
  &=\int_\Omega\boldsymbol{\sigma}^-:\mathbf e(\mathbf v)\dd\mathbf x
    +\int_\Omega q_\varepsilon^+[\mathbf{u}]\,
      \Aeps\bigl[\boldsymbol{\sigma}^+(\mathbf{u}):
      \mathbf{e}(\mathbf{v})\bigr]\dd\mathbf{x}\\
  &=\int_\Omega\bigl(\boldsymbol{\sigma}^-+
      (\Aepsstar q_\varepsilon^+)\boldsymbol{\sigma}^+\bigr)
      :\mathbf{e}(\mathbf{v})\dd\mathbf{y},
\end{align}
which proves \eqref{eq:first-variation}.
\end{proof}

For the exponential potential, we have
\begin{equation}
  q_\varepsilon^+=\exp\left(-\frac{2s_\varepsilon^+}{G_c}\right),
  \quad
  a_\varepsilon^+=\Aepsstar\left[
     \exp\left(-\frac{2\varepsilon}{G_c}
        \Aeps W^+(\mathbf e(\mathbf u))\right)\right].
  \label{eq:exp-coefficient}
\end{equation}

\subsection{Elastodynamic equations and energy balance}
\noindent Let $\rho(\mathbf x)>0$ denote the reference mass density.
The kinetic energy and action functional are defined by
\begin{equation}
 \mathcal T[\dot{\mathbf u}]
   =\frac12\int_\Omega\rho|\dot{\mathbf u}|^2\dd\mathbf x,
 \quad
 \mathcal S[\mathbf u]
   =\int_0^T\left(\mathcal T[\dot{\mathbf u}]
          -\Pi_\varepsilon^{\mathrm{split}}[\mathbf u;t]\right)\dd t.
 \label{eq:dynamic-action}
\end{equation}
Stationarity of $\mathcal S$ for variations vanishing on $\Gamma_D$ and
at the temporal endpoints, followed by integration by parts in time,
gives
\begin{equation}
  \int_\Omega\rho\ddot{\mathbf u}\cdot\mathbf v\dd\mathbf x
  +\int_\Omega\boldsymbol\sigma_\varepsilon^{\mathrm{split}}
      :\mathbf e(\mathbf v)\dd\mathbf x
  =\ell_{\mathrm{ext}}(t)[\mathbf v]
  \quad\text{for all }\mathbf v=\mathbf0\text{ on }\Gamma_D.
  \label{eq:weak-equilibrium}
\end{equation}

\begin{remark}[Crack indicator]
\label{rem:crack-indicator}

We use $1-q_\varepsilon^+[\mathbf u]$ as a crack indicator. For the
exponential potential,
\begin{equation}
 1-q_\varepsilon^+[\mathbf u]
 =1-\exp\!\left(-\frac{2s_\varepsilon^+[\mathbf u]}{G_c}\right)
 \in[0,1).
 \label{eq:crack-indicator}
\end{equation}
It vanishes at zero averaged tensile energy and approaches one in regions
where $q_\varepsilon^+$ tends to zero. The indicator is computed from
displacement and introduces no additional unknown. 
\end{remark}

For sufficiently regular solutions, \eqref{eq:weak-equilibrium}
corresponds to the momentum balance
\begin{equation}
  \rho\ddot{\mathbf u}
    =\operatorname{div}\boldsymbol\sigma_\varepsilon^{\mathrm{split}}
    +\mathbf b
  \quad\text{in }\Omega\times(0,T),
  \label{eq:complete-dynamic-system}
\end{equation}
subject to the prescribed displacement and traction data, with
initial conditions $\mathbf u(\cdot,0)=\mathbf u_0$ and
$\dot{\mathbf u}(\cdot,0)=\mathbf v_0$ compatible with the Dirichlet data.

\begin{remark}[Time-dependent Dirichlet data]
\label{rem:continuous-moving-dirichlet}
For sufficiently regular solutions with prescribed displacement
$\mathbf u=\bar{\mathbf u}(t)$ on $\Gamma_D$, the total mechanical energy
\begin{equation}
 \E(t)=\mathcal T[\dot{\mathbf u}(t)]+\Fsplit[\mathbf u(t)]
 \label{eq:physical-energy-continuous}
\end{equation}
satisfies
\begin{equation}
 \frac{\mathrm d\E}{\mathrm dt}=\ell_{\mathrm{ext}}(t)[\dot{\mathbf u}]+\mathcal P_D(t),
 \label{eq:continuous-work}
\end{equation}
where $ \mathcal P_D(t)=\int_{\Gamma_D}
(\boldsymbol\sigma_\varepsilon^{\mathrm{split}}\mathbf n)
\cdot\dot{\bar{\mathbf u}}\,\mathrm da$ is the power of the reactions associated with
the prescribed motion, and $\mathbf n$ is the outward unit normal.
\end{remark}

If the Dirichlet data are time-independent, $\mathcal P_D=0$; then $\E$ is conserved in the absence of external
power, and $\mathcal H=\E-\ell_{\mathrm{ext}}[\mathbf u]$ is conserved
for a time-independent load functional. These conservation properties motivate the SAV time discretization
in Section~\ref{sec:sav}.

\section{Sharp-crack consistency and emergent tensile stress degradation}
\label{sec:mechanical-meaning}

\noindent \revision{The potential $\Phi$ determines the opening-crack energy and
tensile stress degradation. We examine the elastic limit, planar
displacement transitions and homogeneous tensile response of the split
model.}

\subsection{Elastic limit and displacement jumps}
\label{sec:gamma}
\noindent For a fixed smooth displacement with bounded strain,
$\Phi(s)=s+o(s)$ as $s\downarrow0$ gives
\begin{equation}
 W^-(\mathbf e(\mathbf u))+
 \frac1\varepsilon\Phi\bigl(\varepsilon\Aeps W^+(\mathbf e(\mathbf u))\bigr)
 \longrightarrow W(\mathbf e(\mathbf u))
 \quad(\varepsilon\downarrow0).
 \label{eq:split-elastic-limit}
\end{equation}
Thus, the classical elastic energy density is recovered as $\varepsilon\to0$.

To examine displacement jumps, consider the planar profile
$\mathbf u_\delta(\mathbf x)=\mathbf a\zeta((\mathbf x\cdot\mathbf n)/\delta)$.
Here $\mathbf a$ is the displacement jump across the plane
$\mathbf x\cdot\mathbf n=0$, $\mathbf n$ is its unit normal, and
$\delta>0$ is the transition thickness. The profile $\zeta$ is smooth and
nondecreasing, with $\zeta(r)=0$ for $r\leq-1/2$ and $\zeta(r)=1$ for
$r\geq1/2$.
The matrix $\mathbf a\odot\mathbf n
=(\mathbf a\otimes\mathbf n+\mathbf n\otimes\mathbf a)/2$ has eigenvalues
$(\mathbf a\cdot\mathbf n\pm|\mathbf a|)/2$, with any remaining
eigenvalues zero. Consequently,
\begin{equation}
 W^-(\mathbf a\odot\mathbf n)
 =\frac{\lambda}{2}\langle\mathbf a\cdot\mathbf n\rangle_-^2
  +\frac{\mu}{4}(|\mathbf a|-\mathbf a\cdot\mathbf n)^2,
 \label{eq:split-layer-density}
\end{equation}
and the compressive energy per unit interface area is
\begin{equation}
 \mathcal E_\delta^-=\frac{C_\zeta}{\delta}W^-(\mathbf a\odot\mathbf n),
 \quad C_\zeta=\int_\R|\zeta'(r)|^2\dd r.
 \label{eq:split-compressive-layer}
\end{equation}
For a fixed opening jump $\mathbf a=\gamma\mathbf n$, $\gamma>0$,
\eqref{eq:split-compressive-layer} gives $\mathcal E_\delta^-=0$, so only
the nonlocal tensile term contributes. Let the transition thickness
$\delta=\delta_\varepsilon$ satisfy $\delta_\varepsilon/\varepsilon\to0$
as $\varepsilon\to0$. For an interior planar patch, we use the infinite
planar extension and full averaging balls. With $t=\mathbf x\cdot\mathbf n$
the normal coordinate, the energy per unit interface area is
\begin{equation}
 \begin{aligned}
 \mathcal E_{\varepsilon,\delta_\varepsilon}^{\mathrm{open}}
 &:=\frac1\varepsilon\int_\R
 \Phi\!\left(\varepsilon
 \Aeps[W^+(\mathbf e(\mathbf u_{\delta_\varepsilon}))](t\mathbf n)
 \right)\dd t\\
 &\longrightarrow 2\Phi(\infty)=G_c.
 \end{aligned}
 \label{eq:split-opening-limit}
\end{equation}
The limit follows because averaging produces a layer of thickness
$2\varepsilon+o(\varepsilon)$, within whose interior $\Phi$ approaches
its upper limit. In contrast, tangential sliding or normal compression
gives $W^-(\mathbf a\odot\mathbf n)>0$, so
$\mathcal E_\delta^-\to\infty$ as $\delta\downarrow0$. \revision{\ref{app:planar-transitions} gives these derivations.}

\begin{remark}
\citet{ScillaSolombrino2021} established a Griffith $\Gamma$-limit for the
unsplit reference energy
\begin{equation}
 \mathcal F_\varepsilon[\mathbf u]
 =\frac1\varepsilon\int_\Omega
 \Phi\!\left(\varepsilon\Aeps[W(\mathbf e(\mathbf u))]\right)\dd\mathbf x.
 \label{eq:general-functional}
\end{equation}
This result does not transfer directly to the
split model: $W^+$ is not coercive, and the local term $W^-$ changes the
energy of transitions involving sliding or compression. A full
$\Gamma$-convergence analysis of the split functional is beyond the scope
of this work.
\end{remark}

\subsection{Emergent constitutive degradation and tensile strength}
\noindent The constitutive stress in \eqref{eq:split-stress} determines the
tensile stress coefficient through
$q(s)=\Phi'(s)$ and $a_\varepsilon^+=\Aepsstar[q(s_\varepsilon^+)]$;
no independent degradation law is prescribed
\citep{LussardiNegri2007}. If $\Phi\in C^1([0,\infty))$ is also concave
and satisfies \eqref{eq:saturation-normalization}, then $q$ is nonnegative
and nonincreasing, with
\begin{equation}
 q(0)=1,\qquad \lim_{s\to\infty}q(s)=0.
 \label{eq:q-properties}
\end{equation}
Moreover, $G_c=2\int_0^\infty q(r)\dd r$. \revision{The potential therefore determines the opening-crack energy and,
through the adjoint average of its derivative, the tensile stress
degradation.}

\revision{For a fixed averaging radius, we determine the tensile strength from
a homogeneous scalar tensile mode.}
\begin{proposition}[Homogeneous scalar tensile response]
\label{prop:strength-general}
Let $W^+(\eta)=W(\eta)=E\eta^2/2$ and $W^-=0$ for $\eta\geq0$,
where $E>0$ is the stiffness of this mode. Let $\varepsilon>0$ and
$\Phi\in C^1([0,\infty))\cap C^2((0,\infty))$ satisfy
\eqref{eq:saturation-normalization}. With $q=\Phi'$, the stress conjugate
to the homogeneous strain is
\begin{equation}
 T(\eta)=E\eta\,q(s),\quad s=\frac{\varepsilon E\eta^2}{2}.
 \label{eq:homogeneous-stress}
\end{equation}
Any interior stationary point $s_c>0$ satisfies
\begin{equation}
 \Phi'(s_c)+2s_c\Phi''(s_c)=0.
 \label{eq:peak-condition}
\end{equation}
For a peak at $s_c$, setting $z_c=s_c/G_c$ gives
\begin{equation}
 \eta_c=\sqrt{\frac{2G_c z_c}{\varepsilon E}},\quad
 T_c=C_\Phi\sqrt{\frac{E G_c}{\varepsilon}},
 \label{eq:strength-coefficient}
\end{equation}
where $ C_\Phi=\sqrt{2z_c}\,\Phi'(G_c z_c)$.
\end{proposition}
\begin{proof}
For a homogeneous field, $\Aeps W=W$, so the energy per unit volume is
$\psi_\varepsilon(\eta)=\varepsilon^{-1}\Phi(\varepsilon E\eta^2/2)$.
Differentiation gives \eqref{eq:homogeneous-stress} and
\begin{equation*}
 \frac{\dd T}{\dd\eta}
 =E\bigl(\Phi'(s)+2s\Phi''(s)\bigr).
\end{equation*}
Setting this derivative to zero gives \eqref{eq:peak-condition}; using
$s_c=\varepsilon E\eta_c^2/2$ then gives \eqref{eq:strength-coefficient}.
\end{proof}

For the exponential potential \eqref{eq:exponential-potential}, the peak
occurs at $s_c=G_c/4$, giving
\begin{equation}
 \eta_c=\sqrt{\frac{G_c}{2\varepsilon E}},\qquad
 T_c=\sqrt{\frac{E G_c}{2\mathrm e\,\varepsilon}}.
 \label{eq:exp-strength}
\end{equation}
The coefficients for other potentials are listed in
Table~\ref{tab:strength-laws}. The scaling
$T_c\propto\sqrt{E G_c/\varepsilon}$ shows that the averaging radius
also sets the homogeneous strength scale.

\section{Structure-preserving time discretization}
\label{sec:sav}

\noindent We apply the SAV method to the entire stored-energy functional
\eqref{eq:split-functional}, including its tensile and compressive
contributions. The formulation uses the exact variational derivative of
this functional and introduces a single scalar auxiliary variable.

\subsection{Equivalent SAV formulation}
\noindent Let $V=\{\mathbf w\in H^1(\Omega;\R^d):
\mathbf w=\mathbf0\text{ on }\Gamma_D\}$ and define
\begin{equation}
  m(\mathbf v,\mathbf w)
  =\int_\Omega\rho\,\mathbf v\cdot\mathbf w\,\mathrm d\mathbf x,
  \quad \norm{\mathbf v}_\rho^2=m(\mathbf v,\mathbf v).
  \label{eq:mass-form}
\end{equation}
We first consider time-independent Dirichlet data; prescribed motion
is addressed below.
The internal-force functional is defined by the first variation of the stored energy:
\begin{equation}
\begin{aligned}
  \pair{\Fu(\mathbf u)}{\mathbf w}
  &:=D\Fsplit[\mathbf u](\mathbf w)\\
  &=\int_\Omega\bigl(\boldsymbol\sigma^-(\mathbf u)
       +a_\varepsilon^+[\mathbf u]\boldsymbol\sigma^+(\mathbf u)\bigr)
       :\mathbf e(\mathbf w)\,\mathrm d\mathbf x.
\end{aligned}
  \label{eq:Fu-op}
\end{equation}
Thus \eqref{eq:weak-equilibrium} becomes
\begin{equation}
  \dot{\mathbf u}=\mathbf v,\quad
  m(\dot{\mathbf v},\mathbf w)
  +\pair{\Fu(\mathbf u)}{\mathbf w}
  =\ell_{\mathrm{ext}}(t)[\mathbf w],\quad \mathbf w\in V.
  \label{eq:sav-weak-system}
\end{equation}

Following the SAV construction
\revision{\citep{ShenXuYang2018,ShenXuYang2019,ChengLiuShen2021}}, set
\begin{equation}
  A(\mathbf u)=\Fsplit[\mathbf u]+C_0,\quad
  r(t)=\sqrt{A(\mathbf u(t))},\quad C_0>0.
  \label{eq:continuous-A}
\end{equation}
Here $C_0$ has units of energy, and $A\geq C_0$ because
$\Fsplit\geq0$. The chain rule gives the equivalent system
\begin{equation}
\begin{cases}
  \dot{\mathbf u}=\mathbf v,\\[2pt]
  m(\dot{\mathbf v},\mathbf w)
   +\dfrac{r}{\sqrt{A(\mathbf u)}}
      \pair{\Fu(\mathbf u)}{\mathbf w}
   =\ell_{\mathrm{ext}}(t)[\mathbf w],\quad \mathbf w\in V,\\[5pt]
  \dot r=\dfrac{1}{2\sqrt{A(\mathbf u)}}
                  \pair{\Fu(\mathbf u)}{\mathbf v}.
\end{cases}
  \label{eq:continuous-sav}
\end{equation}
With $r(0)=\sqrt{A(\mathbf u_0)}$, the last equation implies
$\mathrm d(r-\sqrt{A(\mathbf u)})/\mathrm dt=0$.
Hence $r/\sqrt{A(\mathbf u)}=1$ at the continuous level. Define
\begin{equation}
  \Et(t)=\tfrac12\norm{\mathbf v(t)}_\rho^2+r(t)^2-C_0.
  \label{eq:continuous-modenergy}
\end{equation}
Testing the momentum equation with $\mathbf v$ gives $ \frac{\mathrm d\Et}{\mathrm dt}
=\ell_{\mathrm{ext}}(t)[\mathbf v]$. For compatible initial data, the modified
energy $\Et$ equals the physical energy $\E$ in the continuous formulation.
\revision{This equality is generally not preserved by SAV time discretization
\citep{ChengLiuShen2020}.}

\subsection{Structure-preserving SAV scheme}
\noindent In the discrete formulas below, $A$ denotes
$A_h(\mathbf u_h)=\Fsplith[\mathbf u_h]+C_0$, and $\Fu$ denotes its
exact variational derivative. Let $\tau>0$, $t_n=n\tau$, and write
\begin{equation}
  D_tz^{n+1}=\frac{z^{n+1}-z^n}{\tau},\quad
  \delta z^n=z^{n+1}-z^n,\quad
  z^{n+1/2}=\frac{z^{n+1}+z^n}{2}.
  \label{eq:time-notation}
\end{equation}
For the time-independent Dirichlet data considered here,
the displacement belongs to the affine space with prescribed trace,
and the velocity and test functions have zero Dirichlet trace.
Prescribed motion is treated in
Remark~\ref{rem:sav-moving-dirichlet}.
Initialize $r^0=\sqrt{A(\mathbf u^0)}$ and define
\begin{equation}
  \Et^n=\tfrac12\norm{\mathbf v^n}_\rho^2+(r^n)^2-C_0.
  \label{eq:time-energy}
\end{equation}

For $n\geq1$, extrapolate the displacement to the midpoint $  \mathbf u_*^{n+1/2}=\tfrac32\mathbf u^n-\tfrac12\mathbf u^{n-1}$, and define
\begin{equation}
  A_*^{n+1/2}=A(\mathbf u_*^{n+1/2}),\quad
  \Fustar^{n+1/2}=\Fu(\mathbf u_*^{n+1/2}).
  \label{eq:sav-data}
\end{equation}
With
\begin{equation}
  \xi^{n+1/2}=\frac{r^{n+1/2}}{\sqrt{A_*^{n+1/2}}},
  \quad
  \ell_{\mathrm{ext}}^{n+1/2}
    =\ell_{\mathrm{ext}}(t_n+\tau/2),
  \label{eq:sav-ratio}
\end{equation}
the update is
\begin{equation}
\begin{cases}
  D_t\mathbf u^{n+1}=\mathbf v^{n+1/2},\\[2pt]
  m(D_t\mathbf v^{n+1},\mathbf w)
   +\xi^{n+1/2}\pair{\Fustar^{n+1/2}}{\mathbf w}
   =\ell_{\mathrm{ext}}^{n+1/2}[\mathbf w],\\[5pt]
  D_tr^{n+1}=\dfrac{1}{2\sqrt{A_*^{n+1/2}}}
      \pair{\Fustar^{n+1/2}}{D_t\mathbf u^{n+1}}.
\end{cases}
  \label{eq:sav}
\end{equation}
For the first step, use
\begin{equation}
  \mathbf u_*^{1/2}=\mathbf u^0+\tfrac{\tau}{2}\mathbf v^0.
  \label{eq:sav-start}
\end{equation}

\begin{theorem}[Discrete modified-energy balance]
\label{thm:sav-energy}
For time-independent Dirichlet data on a fixed spatial
discretization, \eqref{eq:sav} satisfies, for every $\tau>0$,
\begin{equation}
  \Et^{n+1}-\Et^n
   =\ell_{\mathrm{ext}}^{n+1/2}[\delta\mathbf u^n].
  \label{eq:sav-energy-balance}
\end{equation}
Thus $\Et$ is conserved for zero loading, and
$\widetilde{\mathcal H}^{n+1}=\widetilde{\mathcal H}^n$
for a time-independent load functional, where
$\widetilde{\mathcal H}^n=\Et^n-\ell_{\mathrm{ext}}[\mathbf u^n]$.
\end{theorem}
\begin{proof}
Test the momentum equation with
$\tau\mathbf v^{n+1/2}=\delta\mathbf u^n$ and multiply the scalar
equation by $2\tau r^{n+1/2}$. Apply
$2(a-b,(a+b)/2)=\norm a^2-\norm b^2$ to obtain
\eqref{eq:sav-energy-balance}.
\end{proof}

\begin{remark}[Discrete work of prescribed motion]
\label{rem:sav-moving-dirichlet}
For time-dependent Dirichlet data, the discrete power supplied by the
displacement constraints is
\begin{equation}
 \mathcal P_{D,h}^{n+1/2}
 =(\mathbf R_D^{n+1/2})^{\mathsf T}\mathbf V_D^{n+1/2},
 \label{eq:sav-boundary-power}
\end{equation}
where $\mathbf V_D^{n+1/2}=\delta\bar{\mathbf U}_D^n/\tau$ contains the prescribed boundary velocities
and $\mathbf R_D^{n+1/2}$ is the reaction exerted by the displacement
constraint on the body:
\begin{equation}
 \mathbf R_D^{n+1/2}
 =\left[
   \mathbf M_{\mathrm{all}}
     \frac{\mathbf V^{n+1}-\mathbf V^n}{\tau}
   +\xi^{n+1/2}\nabla\Fsplith(\mathbf U_*^{n+1/2})
   -\mathbf f_{\mathrm{ext}}^{n+1/2}
  \right]_D.
 \label{eq:sav-dirichlet-reaction}
\end{equation}
Here $[\cdot]_D$ selects the prescribed degrees of freedom.
$\mathbf M_{\mathrm{all}}$ is the full consistent mass matrix, and
$\mathbf U$, $\mathbf V$ and $\mathbf f_{\mathrm{ext}}$ contain the
components associated with both free and prescribed degrees of freedom.
The exact discrete balance for the modified energy becomes
\begin{equation}
 \Et^{n+1}-\Et^n
 =\ell_{\mathrm{ext}}^{n+1/2}[\delta\mathbf u^n]
  +\tau\mathcal P_{D,h}^{n+1/2}.
 \label{eq:sav-boundary-work}
\end{equation}
This is the discrete counterpart of \eqref{eq:continuous-work}.
For fixed Dirichlet data, $\mathcal P_{D,h}^{n+1/2}=0$, recovering
Theorem~\ref{thm:sav-energy}.
\end{remark}

With fixed Dirichlet data,
the midpoint velocity system on the free degrees of freedom is
\begin{equation}
\begin{aligned}
 \left(\mathbf M+\frac{\tau^2}{8}\mathbf b_*\mathbf b_*^{\mathsf T}\right)
 \mathbf V^{n+1/2}
 &=\mathbf M\mathbf V^n+\frac{\tau}{2}
       \left(\mathbf f_{\mathrm{ext}}^{n+1/2}-r^n\mathbf b_*\right),\\
 \mathbf b_*&=\frac{\nabla\Fsplith(\mathbf U_*^{n+1/2})}
                  {\sqrt{A_*^{n+1/2}}},
\end{aligned}
 \label{eq:sav-rank-one}
\end{equation}
where $\mathbf M$ is the mass matrix restricted to the free degrees of
freedom, $\mathbf U$ and $\mathbf V$ are the corresponding coefficient
vectors, and the gradient is restricted to the same degrees of freedom.
The new state follows from
\begin{equation}
\begin{aligned}
 \mathbf U^{n+1}&=\mathbf U^n+\tau\mathbf V^{n+1/2},\\
 \mathbf V^{n+1}&=2\mathbf V^{n+1/2}-\mathbf V^n,\\
 r^{n+1}&=r^n+\frac{\tau}{2}\mathbf b_*^{\mathsf T}\mathbf V^{n+1/2}.
\end{aligned}
 \label{eq:sav-recovery}
\end{equation}

\begin{remark}[Modified and physical energies]
The exact balance \eqref{eq:sav-energy-balance} concerns the modified
energy $\Et^n$. The physical discrete energy and the difference between
the two energies are
\begin{equation}
 \begin{aligned}
  \E_h^n&=\tfrac12\norm{\mathbf v^n}_\rho^2+\Fsplith[\mathbf u^n],\\
  \Delta_{\mathrm{SAV}}^n&=\Et^n-\E_h^n
                         =(r^n)^2-A(\mathbf u^n).
 \end{aligned}
 \label{eq:sav-consistency-defect}
\end{equation}
Although initialization gives $\Delta_{\mathrm{SAV}}^0=0$, the scalar
$r^n$ is advanced by its own discrete equation, which does not enforce
$(r^n)^2=A(\mathbf u^n)$ at subsequent steps. Hence this SAV consistency
defect need not remain zero. Under the assumptions of
Theorem~\ref{thm:sav-energy}, the discrete physical energy satisfies
\begin{equation}
 \E_h^{n+1}-\E_h^n
 =\ell_{\mathrm{ext}}^{n+1/2}[\delta\mathbf u^n]
  -\bigl(\Delta_{\mathrm{SAV}}^{n+1}-\Delta_{\mathrm{SAV}}^n\bigr).
 \label{eq:sav-physical-energy-balance}
\end{equation}
Thus the change in $\Delta_{\mathrm{SAV}}^n$ quantifies the discrepancy
between the physical-energy increment and the external work.
\end{remark}

\Needspace{6\baselineskip}
\section{Numerical validation}
\label{sec:numerics}

\noindent We assess the crack-path predictions and propagation dynamics of the
model, together with the SAV scheme of Section~\ref{sec:sav}, using two
benchmarks: dynamic crack branching and Kalthoff--Winkler impact.
For both benchmarks, the explicit velocity-Verlet algorithm
\citep{SwopeEtAl1982}, in kick--drift--kick form, serves as the reference
integrator. Table~\ref{tab:benchmark-materials} summarizes the material
parameters and discretization settings.

Both integrators use adaptive
meshes of bilinear quadrilateral ($Q_1$) elements
\revision{\citep{ShiXuZhang2008}}, with a consistent mass
matrix for SAV and a lumped mass matrix for velocity-Verlet.
The coarsest element size is denoted by $h_0$; two refinement levels yield
$h_{\min}=h_0/4$. The nonlocal averaging radius in each element $K$ is
set to $\varepsilon_K=Mh_K$.

\begin{table}[!htbp]
\centering
\renewcommand{\baselinestretch}{1}
\caption{Material parameters, spatial discretizations and time integrators
for the two dynamic benchmarks.}
\label{tab:benchmark-materials}
\footnotesize
\setlength{\tabcolsep}{4pt}
\renewcommand{\arraystretch}{1.15}
\begin{tabular}{@{}>{\raggedright\arraybackslash}p{0.32\linewidth}
                    >{\raggedright\arraybackslash}p{0.28\linewidth}
                    >{\raggedright\arraybackslash}p{0.32\linewidth}@{}}
\toprule
Quantity & Branching & Kalthoff--Winkler\\
\midrule
$E$ (GPa) & $32$ & $190$\\
$\nu$ & $0.20$ & $0.30$\\
$G_c$ (N/mm) & $0.003$ & $22.13$\\
Thickness $t_p$ (mm) & $1$ & $1$\\
Density $\rho$ (kg/m$^3$) & $2450$ & $8000$\\
Kinematics & Plane strain & Plane strain\\
Elements & Adaptive $Q_1$ & Adaptive $Q_1$\\
Coarse-mesh size $h_0$ (mm)
  & $0.667,\ 0.444,$\newline $0.222$
  & $1.5,\ 1,\ 0.5$\\
Nonlocal radius
  & $\varepsilon_K=5h_K$
  & $\varepsilon_K=3h_K$\\
Integrators & Velocity-Verlet, SAV & Velocity-Verlet, SAV\\
\bottomrule
\end{tabular}
\end{table}

\subsection{Dynamic crack branching}
\label{sec:dynamic-branching}

\begin{figure}[H]
\centering
\scalebox{0.90}{\input{figures/dynamic_branching_setup.tex}}
\caption{Dynamic branching specimen and step loading.}
\label{fig:branching-setup}
\end{figure}
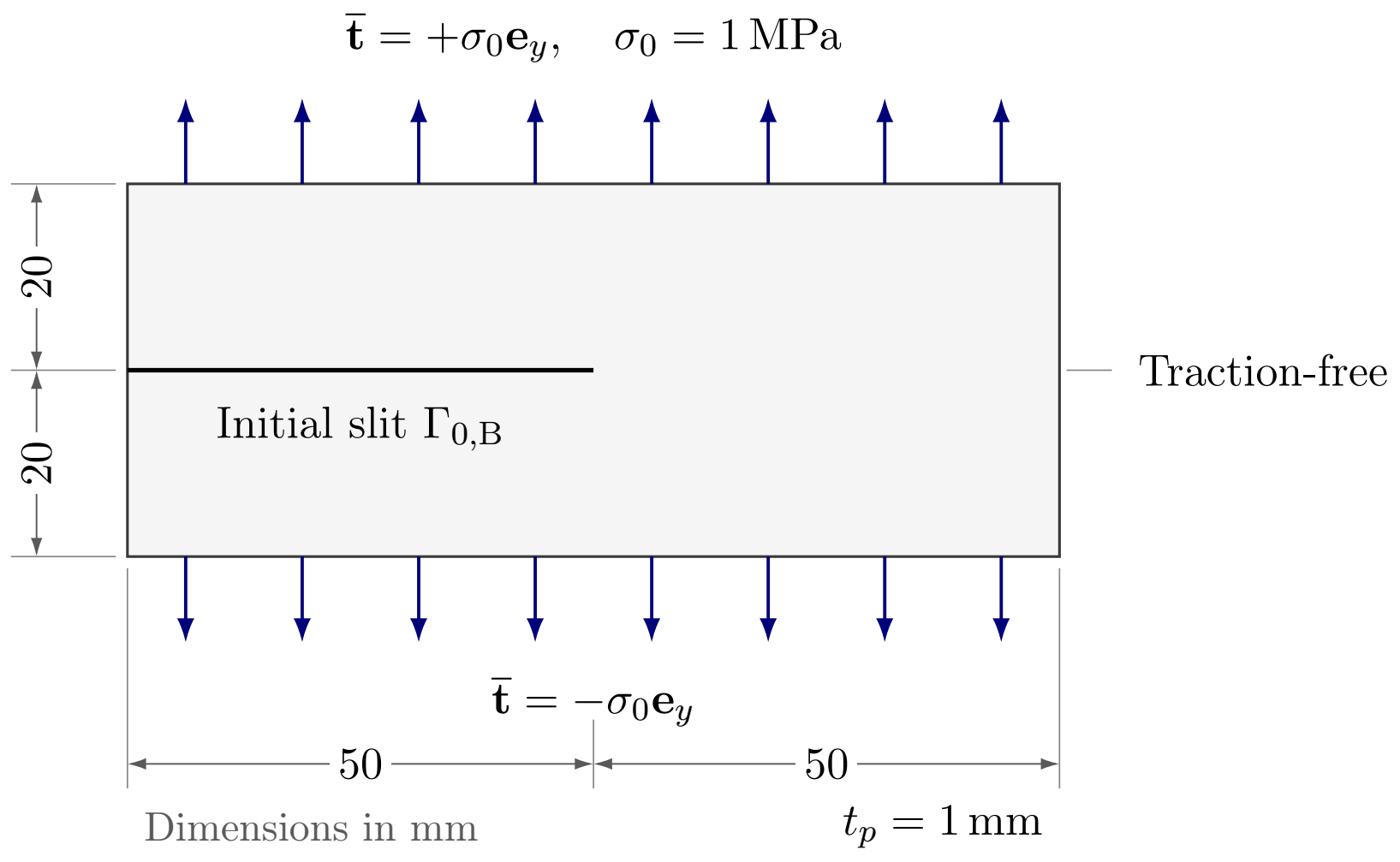

\noindent We simulate dynamic crack branching using the geometry,
material data and loading of \citet{BordenEtAl2012}. With the origin
at the lower-left corner, the plate and initial slit in
Figure~\ref{fig:branching-setup} are
\begin{equation}
\begin{aligned}
 \Omega_{\mathrm B}
  &=\bigl((0,100)\times(0,40)\bigr)\setminus\Gamma_{0,\mathrm B},\\
 \Gamma_{0,\mathrm B}&=[0,50]\times\{20\},
\end{aligned}
\label{eq:branching-geometry}
\end{equation}
with dimensions in millimetres.
At $t=0^+$, equal and opposite $1$ MPa tensile tractions are applied
to the horizontal edges and held constant; the remaining boundaries
are traction-free. Each simulation runs to $80\,\mu\mathrm s$.

\begin{figure}[!htbp]
\centering
\includegraphics[width=\linewidth]{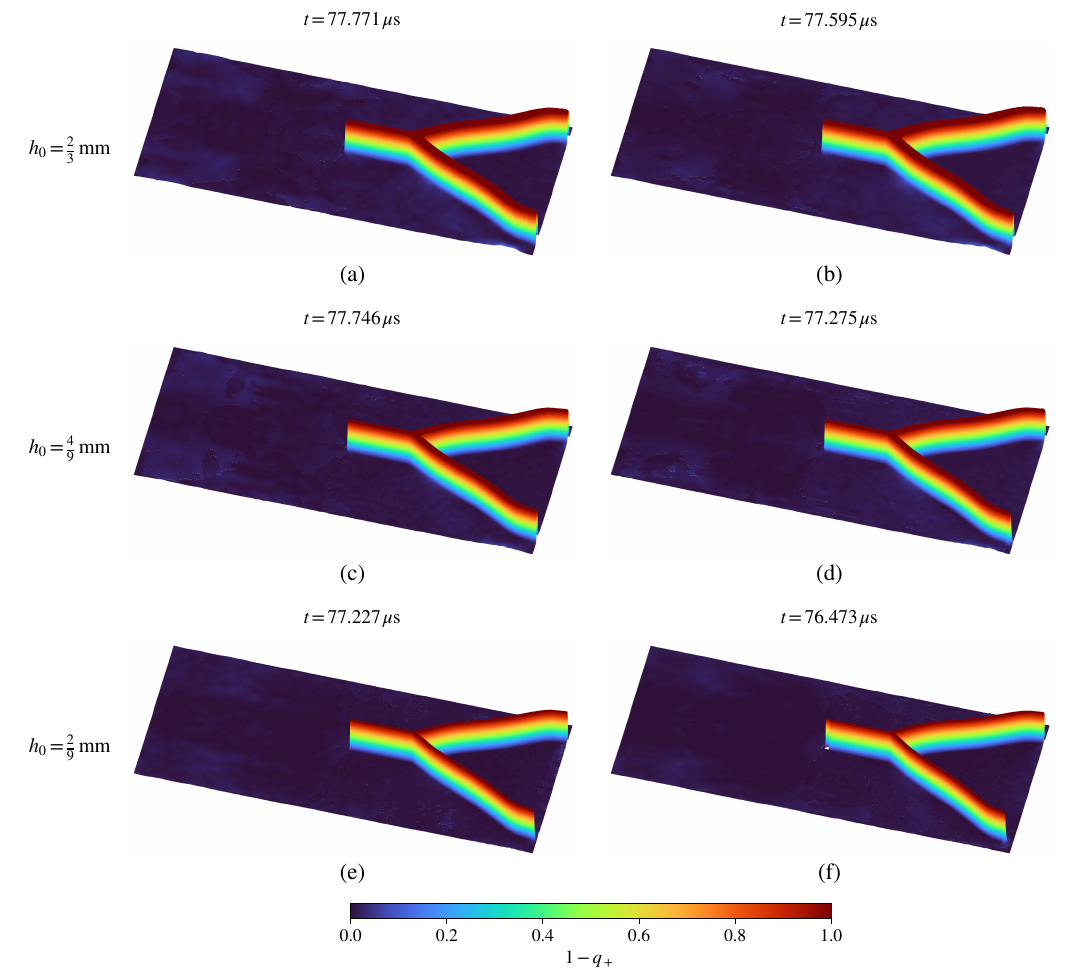}
\caption{Dynamic branching when the connected crack region
$1-q_+\geq0.95$ first reaches the external boundary.
(a,c,e) Velocity-Verlet; (b,d,f) SAV.
The initial mesh size $h_0$ is indicated beside each row. Height and colour both represent
$1-q_+\in[0,1]$ over the reference plane. }
\label{fig:branching-indicator-h0}
\end{figure}

\begin{figure}[!htbp]
\centering
\includegraphics[width=\linewidth]{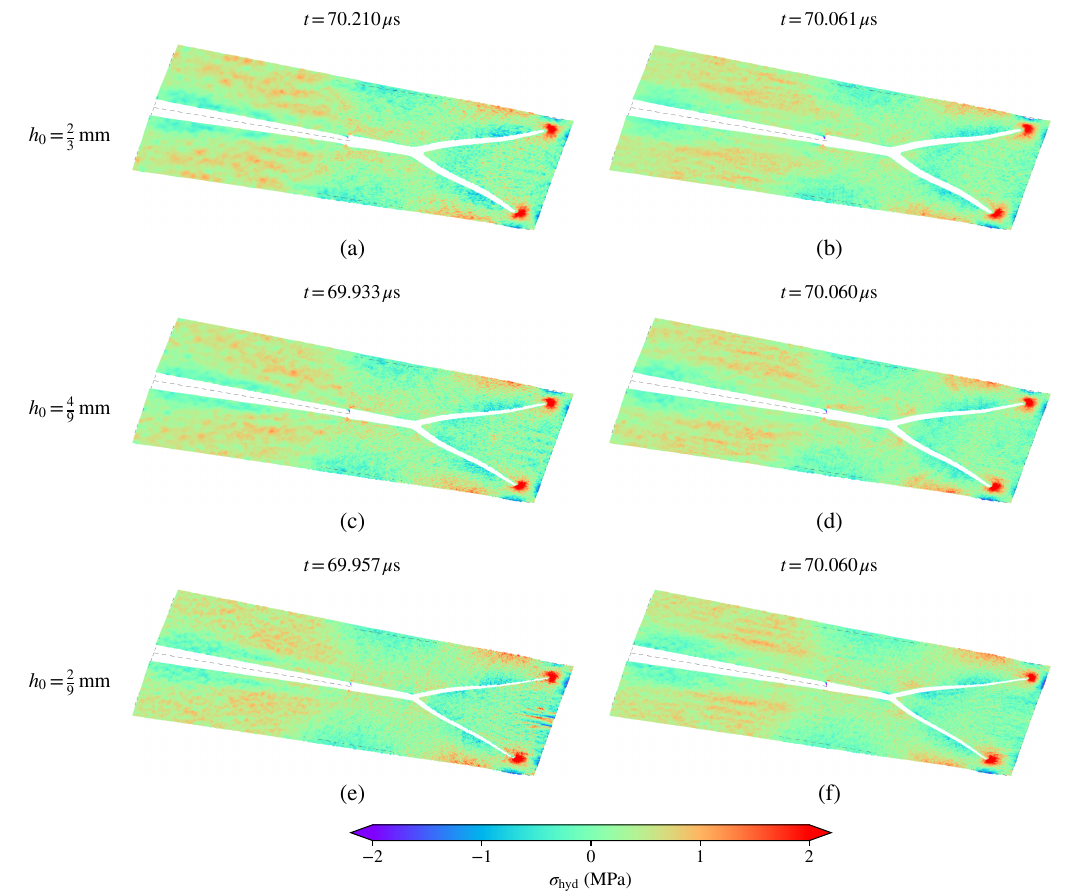}
\caption{Hydrostatic stress
$\sigma_{\mathrm{hyd}}=(\sigma_{xx}+\sigma_{yy}+\sigma_{zz})/3$
near $70\,\mu\mathrm s$. The actual saved times are indicated above each panel.
(a,c,e) Velocity-Verlet; (b,d,f) SAV.
The initial mesh size $h_0$ is indicated beside each row.
The in-plane display coordinates are $\mathbf X+50\mathbf u$.
Height and colour represent the signed stress, with a common height scale
across the six panels. Dashed lines show the initial boundary and notch
in the zero-stress plane. Elements with crack indicator $1-q_+=1$ at all
four nodes are hidden. The common colour scale is $[-2,2]$ MPa.
The ranges of nodal stress values, from top to bottom, are $[-3.95,2.60]$,
$[-5.20,3.33]$ and $[-8.39,4.78]$ MPa in (a,c,e), and
$[-3.74,2.66]$, $[-5.01,3.50]$ and $[-9.83,4.69]$ MPa in (b,d,f).
Values outside the colour scale use its endpoint colours.}
\label{fig:branching-hydrostatic-h0}
\end{figure}

Figure~\ref{fig:branching-indicator-h0} compares the crack patterns
when the connected crack region first reaches the external boundary.
Both integrators consistently reproduce an initially straight extension
followed by two approximately symmetric branches at all three mesh
sizes. The dominant branching pattern is retained as $h_0$ decreases,
while the indicator band narrows. The hydrostatic stress fields near
$70\,\mu\mathrm s$ (Figure~\ref{fig:branching-hydrostatic-h0}) show
tensile stress concentrated near the advancing branch tips.

\begin{table}[H]
\centering
\caption{Dynamic branching: onset, branching times and early
branch angles.}
\label{tab:branching-initiation-angles}
\footnotesize
\setlength{\tabcolsep}{4pt}
\renewcommand{\arraystretch}{1.12}
\begin{tabular}{@{}llcccc@{}}
\toprule
$h_0$ (mm) & Integrator & $t_{\rm det}$ ($\mu\mathrm s$)
 & $t_{\rm br}$ ($\mu\mathrm s$)
 & $\theta_{\rm upper}$ ($^\circ$) & $\theta_{\rm lower}$ ($^\circ$)\\
\midrule
$2/3$ & Velocity-Verlet & $19.103$ & $37.340$ & $27.34$ & $-27.34$\\
      & SAV & $16.579$ & $36.800$ & $28.03$ & $-27.97$\\
\addlinespace
$4/9$ & Velocity-Verlet & $15.566$ & $36.538$ & $28.77$ & $-28.77$\\
      & SAV & $17.500$ & $36.000$ & $29.53$ & $-29.95$\\
\addlinespace
$2/9$ & Velocity-Verlet & $15.141$ & $35.698$ & $30.36$ & $-30.84$\\
      & SAV & $15.143$ & $35.110$ & $30.06$ & $-31.74$\\
\bottomrule
\end{tabular}
\end{table}

\revision{Across the three mesh sizes, the two integrators differ by at most
$0.59\,\mu\mathrm s$ in branching time and $1.2^\circ$ in initial branch
angle (Table~\ref{tab:branching-initiation-angles}).} SAV predicts branch-angle magnitudes
of approximately $28^\circ$--$32^\circ$. The table also lists the first
resolved extension time, $t_{\rm det}$. These measurements use the common
threshold and procedure in \ref{app:crack-tip-speed}.

\FloatBarrier

\begin{figure}[!htbp]
\centering
\includegraphics[width=\linewidth]{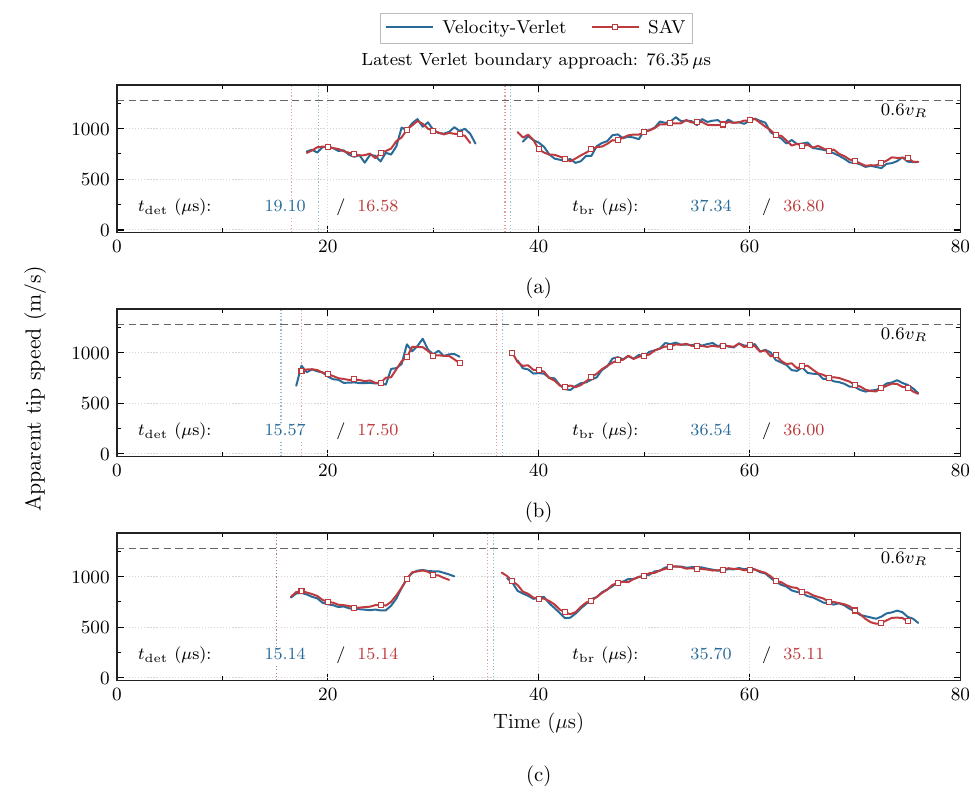}
\caption{Dynamic branching: apparent parent-tip and upper-branch
speeds for velocity-Verlet and SAV at
(a) $h_0=2/3$ mm, (b) $h_0=4/9$ mm, and (c) $h_0=2/9$ mm.}
\label{fig:branching-speed-energy-h0}
\end{figure}

Figure~\ref{fig:branching-speed-energy-h0} shows that the two
integrators also give comparable crack-tip speeds during propagation.
The parent-crack tip speeds obtained with velocity-Verlet reach
approximately $1.0$--$1.2$ km/s. The SAV upper-branch speeds are
approximately $0.96$--$1.10$ km/s over $50$--$60\,\mu\mathrm s$,
close to the corresponding velocity-Verlet speeds, and then decrease
as the tips approach the outer boundary.

\begin{figure}[!htbp]
\centering
\includegraphics[width=\linewidth]{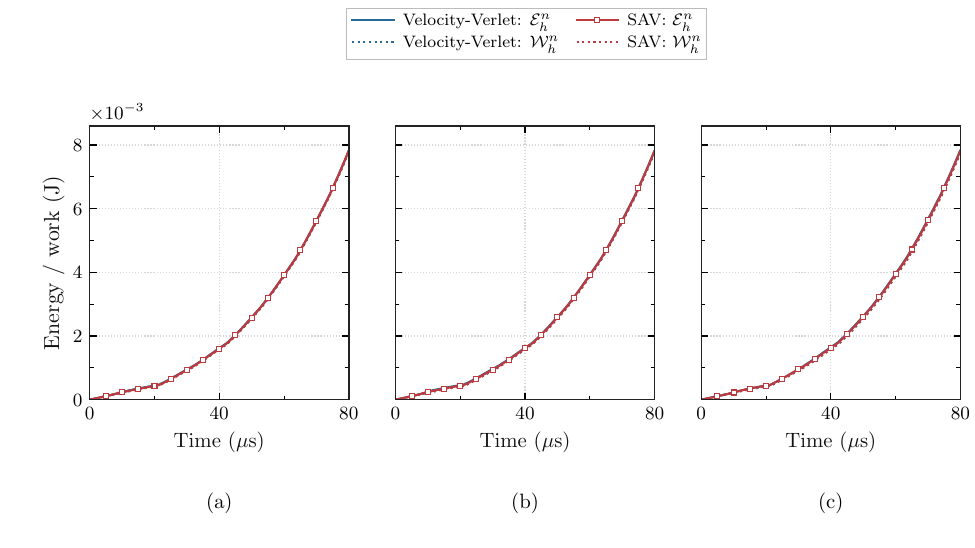}
\caption{Dynamic branching: physical-energy increments $\mathcal E_h^n$ and external
work for velocity-Verlet and SAV at
(a) $h_0=2/3$ mm, (b) $h_0=4/9$ mm, and (c) $h_0=2/9$ mm.}
\label{fig:branching-energy-work}
\end{figure}

\begin{figure}[!htbp]
\centering
\includegraphics[width=0.85\linewidth]{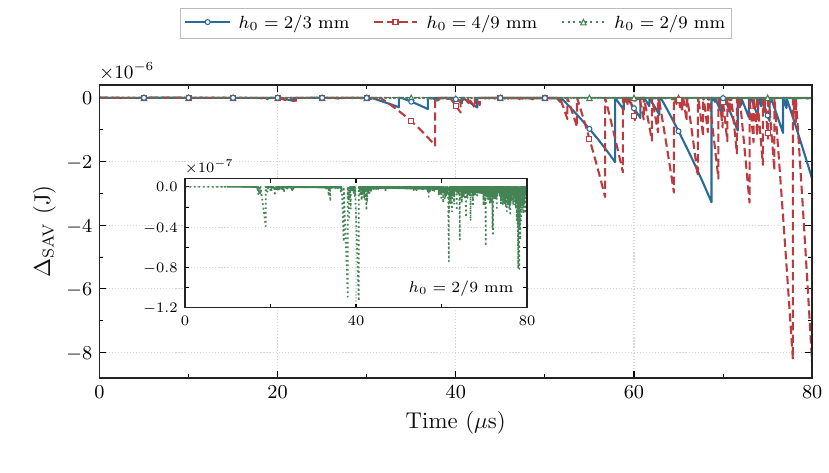}
\caption{SAV consistency defect
$\Delta_{\mathrm{SAV}}^n$ in dynamic branching. The inset enlarges the
$h_0=2/9$ mm curve. The maximum absolute defects
$\max_n|\Delta_{\mathrm{SAV}}^n|$ are $3.28\times10^{-6}$,
$8.35\times10^{-6}$ and $1.13\times10^{-7}$ J for
$h_0=2/3$, $4/9$ and $2/9$ mm, respectively.}
\label{fig:branching-sav-defect}
\end{figure}

\revision{For both integrators and all three mesh sizes, the energy increments
closely follow the external work through branching
(Figure~\ref{fig:branching-energy-work}).} Figure~\ref{fig:branching-sav-defect} reports the SAV consistency
defect; the stepwise modified-energy balance does not require this
defect to vanish.

\FloatBarrier
\subsection{Kalthoff--Winkler impact}
\label{sec:kalthoff-winkler}

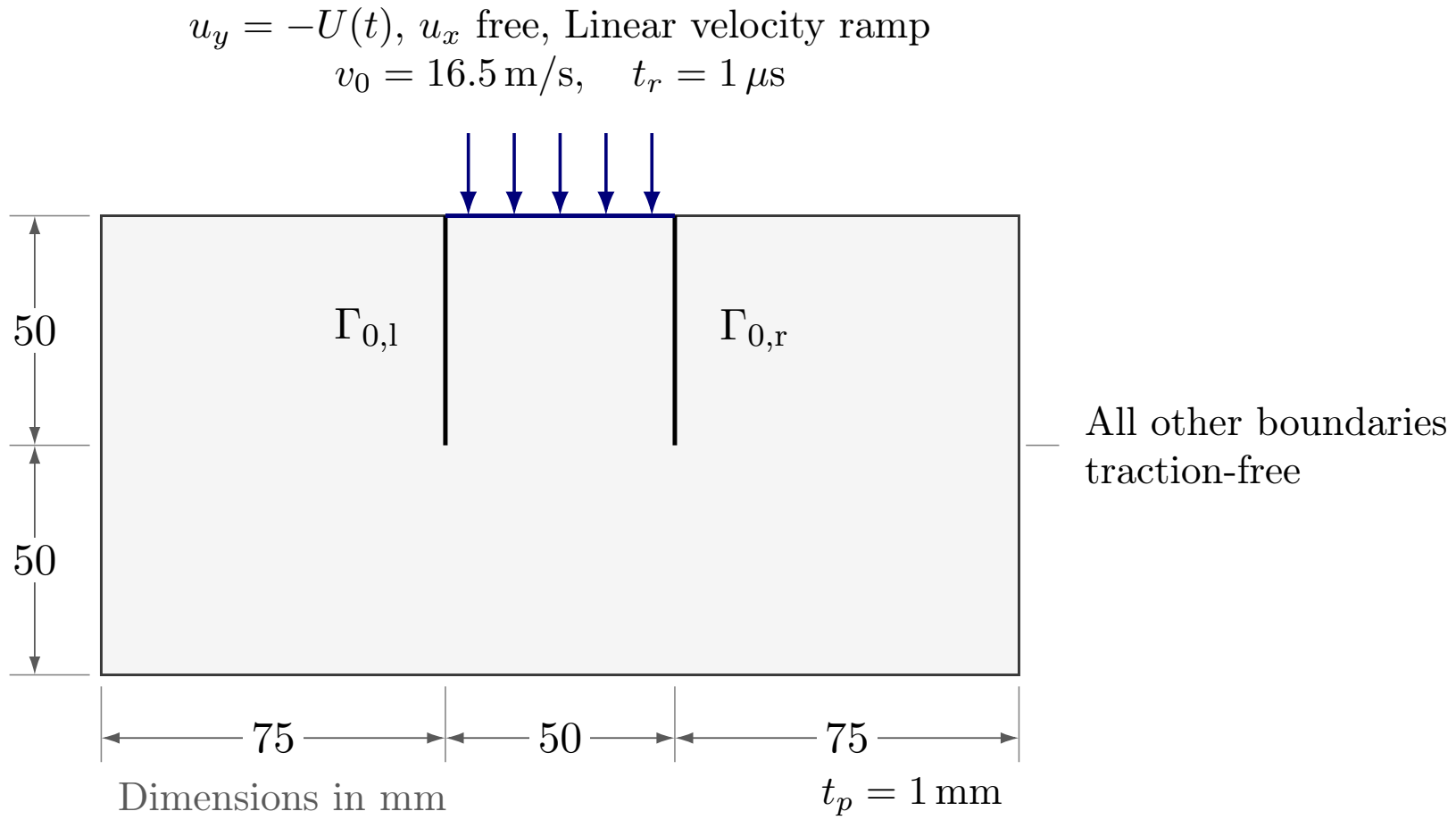
\begin{figure}[!htbp]
\centering
\scalebox{0.85}{\input{figures/kalthoff_winkler_setup.tex}}
\caption{Kalthoff--Winkler geometry and loading. The velocity in \eqref{eq:kw-impact-velocity} is prescribed on the
upper edge between the initial slits; all other boundaries are traction-free. }
\label{fig:kw-setup}
\end{figure}

\noindent This benchmark examines brittle crack kinking under dynamic
shear loading \citep{Kalthoff2000}. The full plate is modelled, with the
loading idealization of \citet{BordenEtAl2012}. With the origin at the
lower-left corner, $x$ horizontal and $y$ vertical, the geometry in
Figure~\ref{fig:kw-setup} is
\begin{equation}
\begin{aligned}
 \Omega_{\mathrm{KW}}
  &=\bigl((0,200)\times(0,100)\bigr)
     \setminus(\Gamma_{0,\mathrm l}\cup\Gamma_{0,\mathrm r}),\\
 \Gamma_{0,\mathrm l}&=\{75\}\times[50,100],\quad
 \Gamma_{0,\mathrm r}=\{125\}\times[50,100],
\end{aligned}
\label{eq:kw-geometry}
\end{equation}
with dimensions in millimetres. A downward velocity is prescribed on
the upper edge of the central region between the initial slits, $\Gamma_v=[75,125]\times\{100\}$:
\begin{equation}
 \dot U(t)=-\dot u_y(t)=v_0
 \begin{cases}
  t/t_r,&0\leq t\leq t_r,\\
  1,&t>t_r,
 \end{cases}
 \quad
 v_0=16.5\ \mathrm{m/s},\quad t_r=1\,\mu\mathrm s.
 \label{eq:kw-impact-velocity}
\end{equation}
Here $U(0)=0$. The horizontal displacement on $\Gamma_v$ is unconstrained; all remaining
boundaries are traction-free. Loading is applied directly to the specimen,
without modelling projectile contact. 

\begin{figure}[!htbp]
\centering
\includegraphics[width=\linewidth]{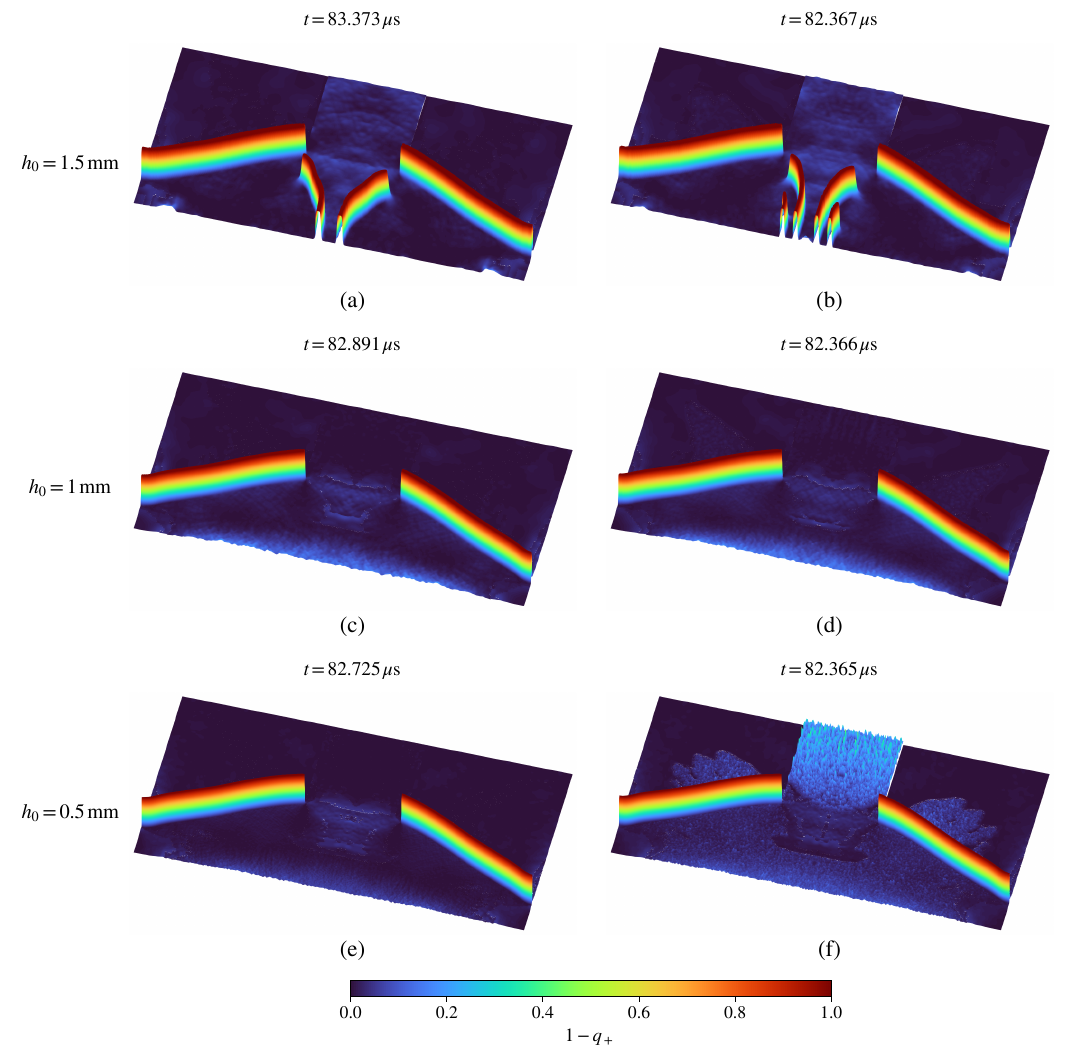}
\caption{Kalthoff--Winkler impact when the connected crack region
$1-q_+\geq0.95$ reaches the external boundary.
(a,c,e) Velocity-Verlet; (b,d,f) SAV.
The initial mesh size $h_0$ is indicated beside each row. Height and colour both represent
$1-q_+\in[0,1]$ over the reference plane; height is a visualization
coordinate.}
\label{fig:kw-indicator-h0}
\end{figure}

\begin{figure}[!htbp]
\centering
\includegraphics[width=\linewidth]{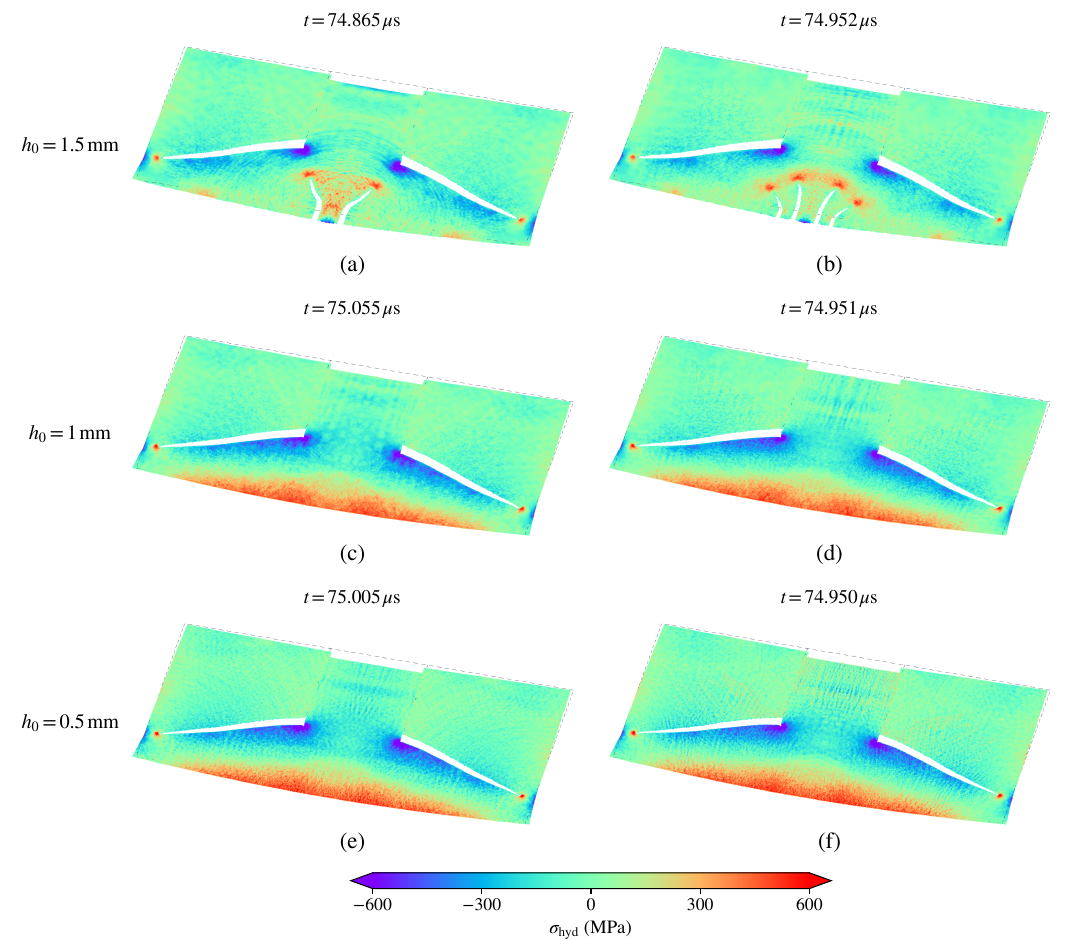}
\caption{Hydrostatic stress
$\sigma_{\mathrm{hyd}}=(\sigma_{xx}+\sigma_{yy}+\sigma_{zz})/3$
near $75\,\mu\mathrm s$. The actual saved times are indicated above each panel.
(a,c,e) Velocity-Verlet; (b,d,f) SAV.
The initial mesh size $h_0$ is indicated beside each row. The in-plane display coordinates are
$\mathbf X+5\mathbf u$. Height and colour represent the signed stress,
with a common height scale across the six panels. Dashed lines show the
initial geometry in the zero-stress plane.
Fully saturated cells are hidden. The common scale is $[-600,600]$ MPa;
nodal ranges from top to bottom are $[-2045,545]$, $[-2339,583]$ and
$[-3284,968]$ MPa in (a,c,e), and $[-2010,604]$, $[-2350,660]$ and
$[-3272,1125]$ MPa in (b,d,f). Out-of-range values use endpoint
colours.}
\label{fig:kw-hydrostatic-h0}
\end{figure}

Figure~\ref{fig:kw-indicator-h0} shows that both integrators reproduce
two approximately symmetric outward-kinking cracks across all three
initial mesh sizes. This dominant pattern persists as $h_0$ decreases,
while the indicator bands narrow. The coarsest cases also show additional
localization at the opposite boundary, $y=0$.
Figure~\ref{fig:kw-hydrostatic-h0} compares the hydrostatic stress fields
near $75\,\mu\mathrm s$. All time-history comparisons below end at
$83\,\mu\mathrm s$, when the cracks are near the outer boundaries.

\FloatBarrier

\begin{table}[!htbp]
\centering
\caption{Kalthoff--Winkler: onset and kink angles.}
\label{tab:kw-initiation-angles}
\small
\setlength{\tabcolsep}{6pt}
\renewcommand{\arraystretch}{1.12}
\begin{tabular}{@{}llccc@{}}
\toprule
$h_0$ (mm) & Integrator & $t_{\rm det}$ ($\mu\mathrm s$)
           & $\theta_{\mathrm l}$ ($^\circ$)
           & $\theta_{\mathrm r}$ ($^\circ$)\\
\midrule
$1.5$ & Velocity-Verlet & $24.241$ & $-74.89$ & $+74.89$\\
      & SAV & $23.997$ & $-74.85$ & $+74.85$\\
\addlinespace
$1$   & Velocity-Verlet & $24.522$ & $-72.52$ & $+72.52$\\
      & SAV & $24.301$ & $-73.19$ & $+73.19$\\
\addlinespace
$0.5$ & Velocity-Verlet & $23.998$ & $-72.65$ & $+72.65$\\
      & SAV & $23.925$ & $-73.45$ & $+73.45$\\
\bottomrule
\end{tabular}
\end{table}

\begin{figure}[!htbp]
\centering
\includegraphics[width=\linewidth]{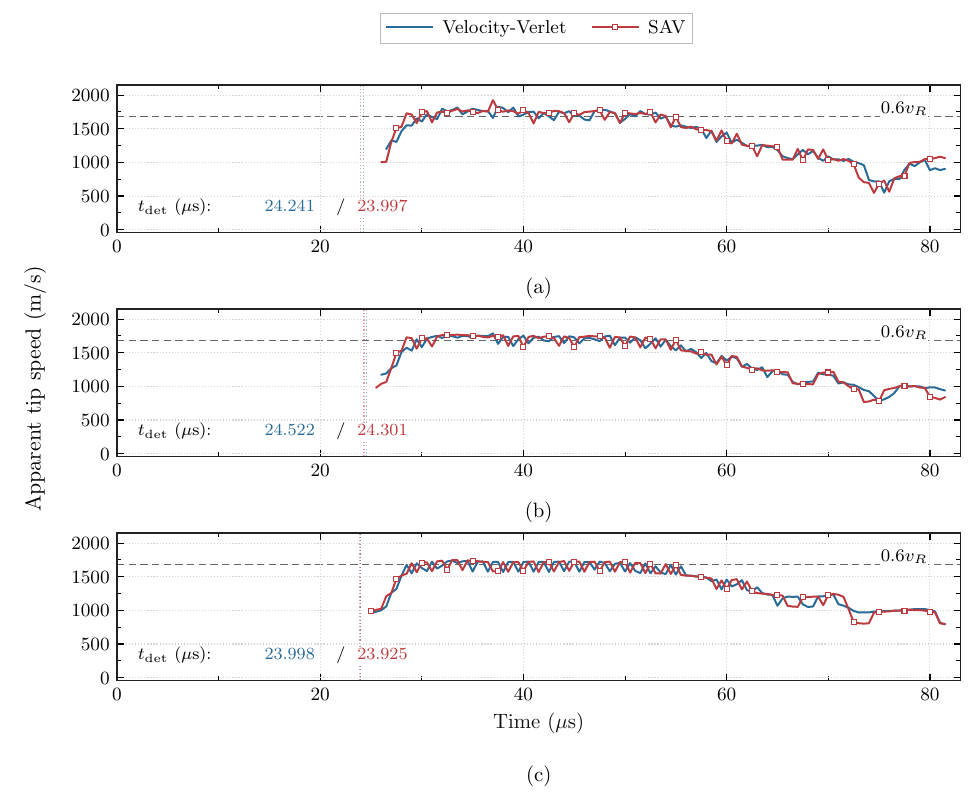}
\caption{Kalthoff--Winkler: apparent right-tip speed up to
$83\,\mu\mathrm s$ for velocity-Verlet and SAV at
(a) $h_0=1.5$ mm, (b) $h_0=1$ mm, and (c) $h_0=0.5$ mm.}
\label{fig:kw-speed-h0}
\end{figure}

\FloatBarrier

Table~\ref{tab:kw-initiation-angles} quantifies the close agreement:
at each resolution, the two integrators differ by less than
$0.25\,\mu\mathrm s$ in the first resolved extension time and at most
$0.8^\circ$ in kink-angle magnitude. SAV predicts extension near
$24\,\mu\mathrm s$ and angle magnitudes of $73.19^\circ$--$74.85^\circ$,
close to the approximately $70^\circ$ angle reported for this benchmark
\citep{BordenEtAl2012}. Negative and positive signs denote clockwise
and counterclockwise kinking from the downward slit direction.
Both integrators use the same measurement procedure
(\ref{app:crack-tip-speed}).

\begin{figure}[!htbp]
\centering
\includegraphics[width=\linewidth]{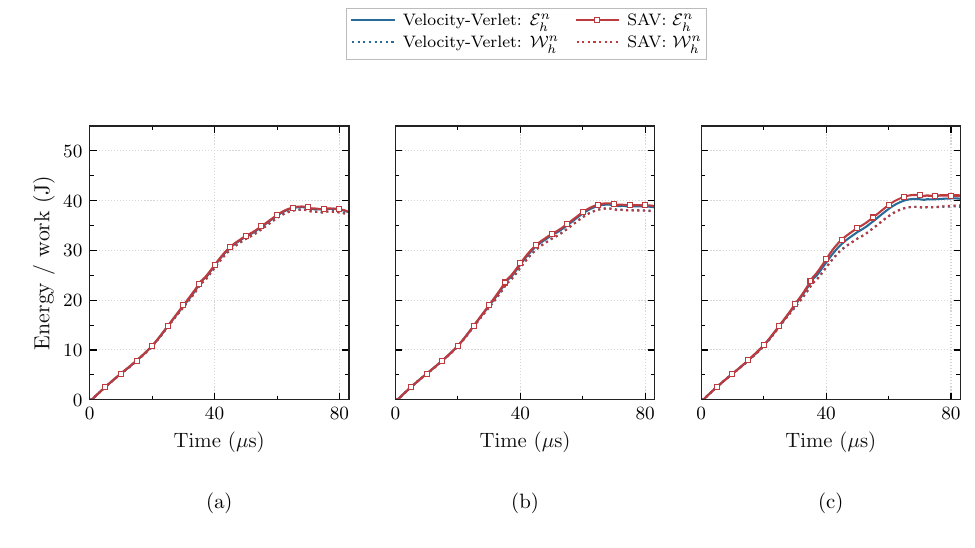}
\caption{Kalthoff--Winkler: physical-energy increments $\mathcal E_h^n$ and external
work for velocity-Verlet and SAV up to $83\,\mu\mathrm s$ at
(a) $h_0=1.5$ mm, (b) $h_0=1$ mm, and (c) $h_0=0.5$ mm.}
\label{fig:kw-energy-work}
\end{figure}

Figure~\ref{fig:kw-speed-h0} shows comparable apparent right-tip
speed histories for velocity-Verlet and SAV at each initial mesh size.
Both integrators capture acceleration to approximately $1.6$--$1.8$ km/s
over $35$--$50\,\mu\mathrm s$, followed by a decrease as the cracks
approach the outer boundaries. \revision{Both integrators give similar histories of energy increments and external
work over the $83\,\mu\mathrm s$ comparison interval
(Figure~\ref{fig:kw-energy-work}).}

\FloatBarrier

We compare the model with tension--compression splitting against the
unsplit model \eqref{eq:general-functional} at $h_0=0.5$ mm, using
otherwise identical physical parameters and numerical settings.

\begin{figure}[!htbp]
\centering
\includegraphics[width=\linewidth]{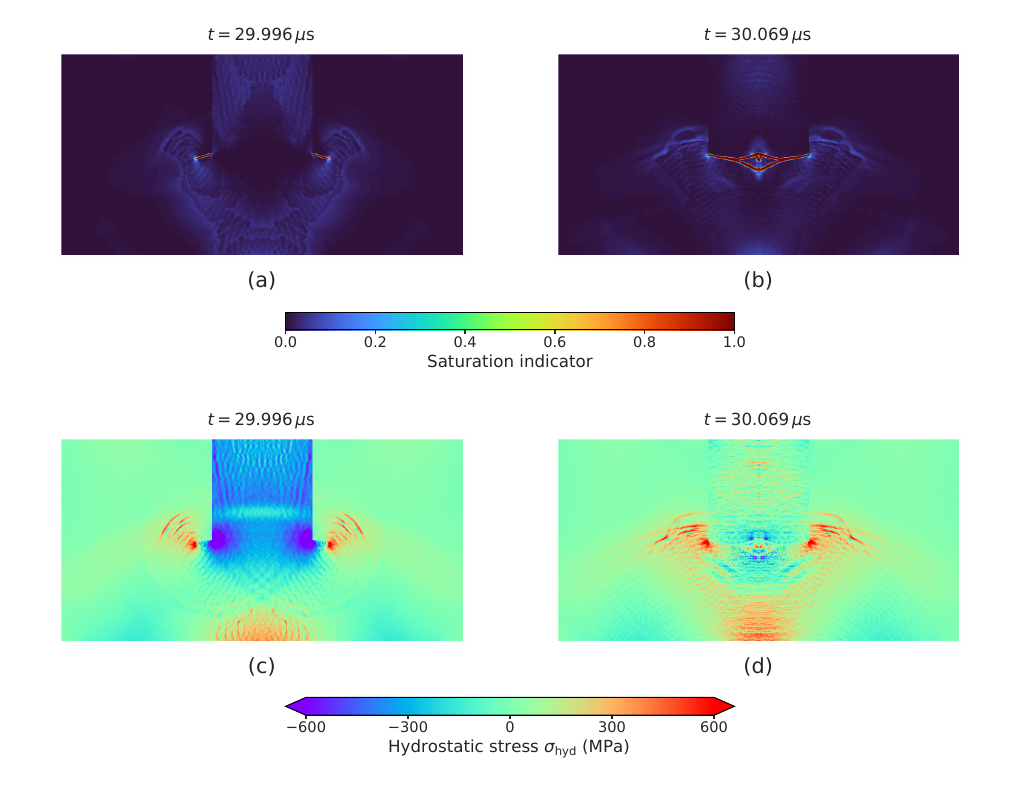}
\caption{Kalthoff--Winkler, $h_0=0.5$ mm: fields computed with velocity-Verlet
for the split and unsplit models near $30\,\mu\mathrm s$ on the reference geometry.
(a,c) Split model; (b,d) unsplit model.
(a,b) Saturation indicators $1-q_+$ and $1-q$, respectively;
(c,d) hydrostatic stress.
All elements are shown. Stress ranges in (c,d) are $[-3559,949]$ and
$[-940,1002]$ MPa, respectively; values outside $[-600,600]$ MPa use endpoint colours.}
\label{fig:kw-split-unsplit-fields}
\end{figure}

Near $30\,\mu\mathrm s$, the split model reproduces outward kinking,
whereas the unsplit model develops an inward localization band between
the notch tips and markedly reduces the compressive stress in the
impacted central region (Figure~\ref{fig:kw-split-unsplit-fields}).

\begin{figure}[!htbp]
\centering
\includegraphics[width=\linewidth]{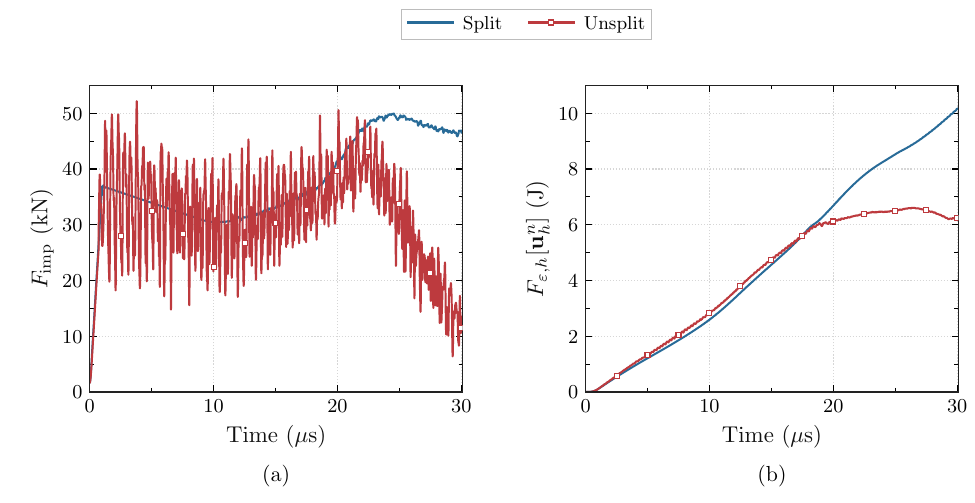}
\caption{Kalthoff--Winkler, $h_0=0.5$ mm: time histories computed with
velocity-Verlet for the split and unsplit models.
(a) Impact reaction $F_{\mathrm{imp}}$;
(b) stored energy $F_{\varepsilon,h}[\mathbf u_h^n]$.}
\label{fig:kw-split-unsplit-history}
\end{figure}

\noindent Figure~\ref{fig:kw-split-unsplit-history} shows the stored energy and
impact reaction
$F_{\mathrm{imp}}^n=\sum_{a\in D_v}
(m_a\ddot u_a^n+f_{\mathrm{int},a}^n)$,
where $D_v$ denotes the prescribed impact degrees of freedom,
$m_a$ the lumped masses, and $f_{\mathrm{int},a}$ the internal forces,
with components taken positive in the impact direction.
\revision{Both the stored energy and impact reaction are substantially smaller
in the unsplit model.}

\FloatBarrier

\revision{The outward crack paths and retained compressive response support
tension--compression splitting for physically plausible mixed-mode
dynamic cracking.}

\FloatBarrier
\Needspace{8\baselineskip}
\section{Conclusions and outlook}
\label{sec:conclusion}

\noindent \revision{We have developed a variational diffuse-crack model based on
nonlocal averaging of tensile strain energy.} A bounded potential acts
on the scaled tensile-energy average, while the compressive energy remains
local. The first variation of the energy determines tensile stress
degradation through the adjoint averaging operator, without an independent
fracture field or a separately prescribed degradation law. At a fixed
averaging radius, the shape of the potential determines the homogeneous
tensile strength.

\revision{The stored-energy functional recovers classical elasticity for fixed
smooth displacement fields and the Griffith surface energy for the planar
pure-opening profiles considered here.}

The linearly implicit SAV scheme satisfies an exact discrete balance
between the modified-energy increment and the external work. The two
dynamic benchmarks illustrate crack branching and kinking. In the
Kalthoff--Winkler benchmark, the unsplit model predicts inward localization
between the notch tips and a smaller compressive reaction. \revision{The split model preserves the compressive elastic response by retaining
the compressive energy locally.} \revision{Further work will address $\Gamma$-convergence of the split functional
and extend the model to multiphysics fracture, fracture in elastoplastic
and hyperelastic solids, and fracture of shells and rods.}

\appendix
\renewcommand{\thetheorem}{\Alph{section}.\arabic{theorem}}

\section{Representative saturation potentials}
\label{app:saturation-laws}

The following saturation potentials satisfy $\Phi'(0)=1$ and
$\Phi(\infty)=G_c/2$. Their derivatives decrease from unity to zero at
rates determined by the form of the potential.

\begin{table}[htbp]
\centering
\caption{\revision{Saturation potentials satisfying
\eqref{eq:saturation-normalization}.}}
\label{tab:saturation-laws}
\renewcommand{\arraystretch}{1.45}
\begin{tabular}{@{}>{\raggedright\arraybackslash}p{0.17\linewidth}
                    >{\raggedright\arraybackslash}p{0.43\linewidth}
                    >{\raggedright\arraybackslash}p{0.30\linewidth}@{}}
\toprule
Family & $\Phi(s)$ & $q(s)=\Phi'(s)$\\
\midrule
Exponential
& $\displaystyle \frac{G_c}{2}\bigl(1-e^{-2s/G_c}\bigr)$
& $\displaystyle e^{-2s/G_c}$\\
Rational
& $\displaystyle \frac{s}{1+2s/G_c}$
& $\displaystyle (1+2s/G_c)^{-2}$\\
Arctangent
& $\displaystyle \frac{G_c}{\pi}\tan^{-1}\!\left(\frac{\pi s}{G_c}\right)$
& $\displaystyle \frac{1}{1+(\pi s/G_c)^2}$\\
Hyperbolic tangent
& $\displaystyle \frac{G_c}{2}\tanh\!\left(\frac{2s}{G_c}\right)$
& $\displaystyle \sech^2\!\left(\frac{2s}{G_c}\right)$\\
\bottomrule
\end{tabular}
\end{table}

\begin{figure}[H]
  \centering
  \includegraphics[width=\textwidth]{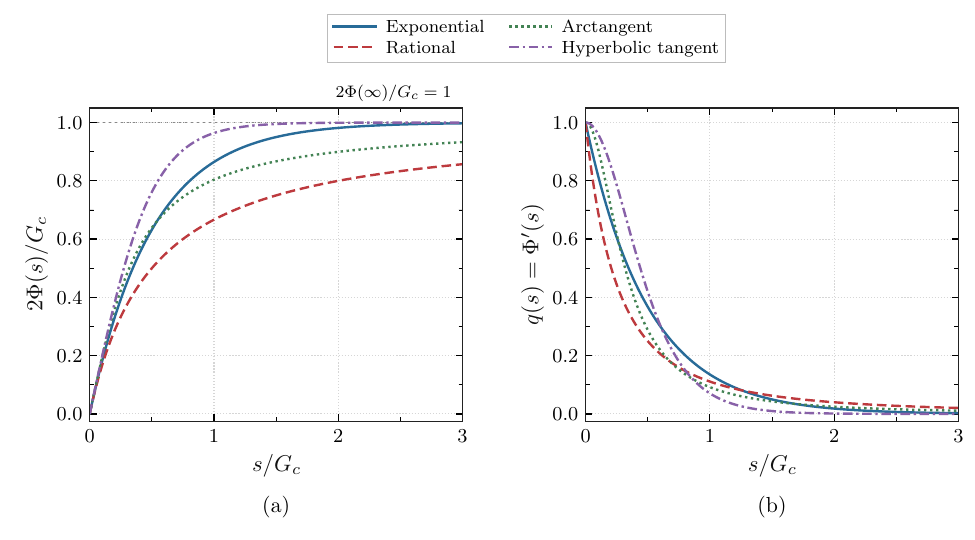}
  \caption{Comparison of (a) the normalized saturation potential
  $2\Phi(s)/G_c$ and (b) the corresponding degradation function
  $q(s)=\Phi'(s)$ for the four potentials in
  \cref{tab:saturation-laws}.}
  \label{fig:saturation-law-comparison}
\end{figure}

Figure~\ref{fig:saturation-law-comparison}
compares the approach of the potentials to their common saturation limit
and the corresponding decay of $q$.

The dimensionless peak locations $z_c=s_c/G_c$ and strength coefficients
$C_\Phi$ in \eqref{eq:strength-coefficient}, obtained from
\cref{eq:peak-condition}, are listed in \cref{tab:strength-laws}.

\begin{table}[H]
\centering
\caption{Homogeneous uniaxial peak locations and strength coefficients
in \eqref{eq:strength-coefficient} for the potentials in
\cref{tab:saturation-laws}.
For the hyperbolic-tangent potential, $z_c$ is the unique positive
root of $8z_c\tanh(2z_c)=1$.}
\label{tab:strength-laws}
\renewcommand{\arraystretch}{1.35}
\begin{tabular}{@{}lll@{}}
\toprule
Family & $z_c=s_c/G_c$ & $C_\Phi$\\
\midrule
Exponential
& $1/4$
& $1/\sqrt{2\mathrm e}\approx0.428882$\\
Rational
& $1/6$
& $9/(16\sqrt{3})\approx0.324760$\\
Arctangent
& $1/(\pi\sqrt{3})$
& $\frac{3}{4}\sqrt{2/(\pi\sqrt{3})}\approx0.454696$\\
Hyperbolic tangent
& $\approx0.260907$
& $\sqrt{2z_c}\,\sech^2(2z_c)\approx0.556558$\\
\bottomrule
\end{tabular}
\end{table}

\FloatBarrier

\section{Planar displacement transitions in the split model}
\label{app:planar-transitions}

\noindent For the profile introduced in Section~\ref{sec:gamma}, set
$t=\mathbf x\cdot\mathbf n$. Differentiation gives
\begin{equation}
 \mathbf e(\mathbf u_\delta)
 =\frac{\zeta'(t/\delta)}{\delta}\mathbf A,
 \quad
 \mathbf A=\mathbf a\odot\mathbf n
 :=\frac{\mathbf a\otimes\mathbf n+\mathbf n\otimes\mathbf a}{2}.
 \label{eq:split-planar-transition}
\end{equation}
Decompose $\mathbf a=a_n\mathbf n+\mathbf a_\tau$, where
$a_n=\mathbf a\cdot\mathbf n$ and $\mathbf a_\tau\cdot\mathbf n=0$.
In a basis aligned with $\mathbf n$ and $\mathbf a_\tau$, the nonzero
part of $\mathbf A$ lies in the block
\begin{equation}
 \mathbf A_{(2)}=
 \begin{pmatrix}a_n&|\mathbf a_\tau|/2\\|\mathbf a_\tau|/2&0\end{pmatrix},
 \quad
 \det(\mathbf A_{(2)}-\eta\mathbf I)
 =\eta^2-a_n\eta-\frac{|\mathbf a_\tau|^2}{4}.
 \label{eq:split-planar-block}
\end{equation}
An arbitrary tangential axis is used when $\mathbf a_\tau=\mathbf0$.
The two eigenvalues are
\begin{equation}
 \eta_\pm=\frac{\mathbf a\cdot\mathbf n\pm|\mathbf a|}{2},
 \qquad \eta_-\leq0\leq\eta_+,
 \label{eq:split-planar-eigenvalues}
\end{equation}
and any remaining eigenvalues vanish. Hence
$\operatorname{tr}\mathbf A=\mathbf a\cdot\mathbf n$ and
$\mathbf A^-:\mathbf A^-=\eta_-^2$, which gives
\eqref{eq:split-layer-density}. Since $\zeta'\geq0$ and $W^-$ is
quadratic under nonnegative scaling, integration along the normal yields
\begin{equation*}
 \int_\R W^-\!\left(\frac{\zeta'(t/\delta)}{\delta}\mathbf A\right)\dd t
 =\frac{W^-(\mathbf A)}{\delta}\int_\R|\zeta'(r)|^2\dd r,
\end{equation*}
which proves \eqref{eq:split-compressive-layer}.
For pure opening, $\mathbf a=\gamma\mathbf n$ with fixed $\gamma>0$,
this term vanishes. For the local calculation per unit interface area,
we use the infinite planar extension with full averaging balls, representing
an interior patch away from the exterior boundary and interface edges.
\revision{For
$\delta_\varepsilon=o(\varepsilon)$, the argument of $\Phi$
diverges at each fixed rescaled normal position $t/\varepsilon\in(-1,1)$.}
The affected layer lies within
$|t|\leq\varepsilon+\delta_\varepsilon/2$; hence boundedness of $\Phi$
and dominated convergence give
\begin{equation}
 \mathcal E_{\varepsilon,\delta_\varepsilon}^{\mathrm{open}}
 =\frac1\varepsilon\int_\R
   \Phi\!\left(\varepsilon\Aeps
       W^+(\mathbf e(\mathbf u_{\delta_\varepsilon}))\right)\dd t
 \longrightarrow\int_{-1}^{1}\Phi(\infty)\dd r=G_c.
 \label{eq:opening-layer-proof}
\end{equation}
Finally, $W^+(-r\mathbf I)=0$ for every $r>0$, which verifies the lack of
coercivity stated in Section~\ref{sec:gamma}.

\section{\revision{Crack-path identification and apparent tip speed}}
\label{app:crack-tip-speed}
\noindent Initiation, branching and angles are extracted from crack-connected Gauss
points satisfying $1-q_+\geq0.95$. The recorded tip $\mathbf r_k$ at time $t_k$ is the ridge tip for
Kalthoff--Winkler and the tip of the connected front of the nodally
projected $Q_1$ indicator for branching, using the same threshold. \revision{In the branching test, the recorded tip follows the parent crack up to
branching and the upper branch thereafter.}
For $t\in[t_k,t_{k+1}]$, define
\begin{equation}
 \widehat{\mathbf r}(t)
 =(1-\alpha)\mathbf r_k+\alpha\mathbf r_{k+1},\qquad
 \alpha=\frac{t-t_k}{t_{k+1}-t_k}.
 \label{eq:app-tip-interpolation}
\end{equation}
The apparent tip speed is
\begin{equation}
 v_{\mathrm{app}}(t;\tau)
 =\frac{\|\widehat{\mathbf r}(t+\tau/2)
              -\widehat{\mathbf r}(t-\tau/2)\|}{\tau},
 \quad \tau=2\,\mu\mathrm s,
 \label{eq:app-apparent-tip-speed}
\end{equation}
evaluated at centres spaced by $0.5\,\mu\mathrm s$.

\section*{Declaration of competing interest}
The authors declare no known competing financial interests or personal
relationships that could have appeared to influence the work reported in
this paper.

\bibliographystyle{elsarticle-num-names}
\begingroup
\setlength{\bibsep}{4pt plus 1pt minus 1pt}
\expandafter\def\csname revisedbib@ChengSunXu2026\endcsname{}
\expandafter\def\csname revisedbib@YuSun2023\endcsname{}
\expandafter\def\csname revisedbib@SunYu2025\endcsname{}
\expandafter\def\csname revisedbib@ChengLiuShen2021\endcsname{}
\expandafter\def\csname revisedbib@ChengLiuShen2020\endcsname{}
\expandafter\def\csname revisedbib@ShiXuZhang2008\endcsname{}
\let\originalbibitem\bibitem
\renewcommand{\bibitem}[2][]{%
  \color{black}\originalbibitem[#1]{#2}%
  \ifcsname revisedbib@#2\endcsname\color{black}\fi}
\color{black}
\bibliography{references}
\endgroup

\end{document}

%% file: figures/split_functional_schematic.tex
\begingroup
\definecolor{sfblue}{RGB}{35,81,116}
\colorlet{sfaccent}{sfblue}
\colorlet{sflight}{white}
\colorlet{sfgray}{black!65}
\begin{tikzpicture}[
  x=1cm,y=1cm,font=\fontsize{9}{11}\selectfont,text=black,
  frame/.style={draw=black!85,line width=.55pt},
  load/.style={draw=black!85,line width=.55pt,-{Stealth[length=1.5mm,width=1.05mm]}},
  axis/.style={draw=black!85,line width=.55pt,-{Stealth[length=1.5mm,width=1.05mm]}},
  note/.style={font=\fontsize{8}{10}\selectfont,text=sfgray}
]
  \path[use as bounding box] (-.1,.05) rectangle (13.35,8.53);

  \foreach \cx in {3.00,10.30} {
    \begin{scope}[shift={(\cx,6.80)}]
      \path[fill=sflight] (0,0) circle[radius=1.52];
      \begin{scope}
        \clip (0,0) circle[radius=1.52];
        \foreach \g in {-1.20,-.80,-.40,0,.40,.80,1.20} {
          \draw[draw=black!12,line width=.18pt]
            ({.88*\g},-1.65)
            .. controls ({1.04*\g},-.55) and ({1.04*\g},.55)
            .. ({.88*\g},1.65);
          \draw[draw=black!12,line width=.18pt]
            (-1.65,{.88*\g})
            .. controls (-.55,{1.04*\g}) and (.55,{1.04*\g})
            .. (1.65,{.88*\g});
        }
      \end{scope}
      \draw[frame] (0,0) circle[radius=1.52];
      \draw[black!85,line width=1.05pt]
        ({1.52*cos(145)},{1.52*sin(145)})
        arc[start angle=145,end angle=215,radius=1.52];
      \draw[draw=sfgray,line width=.45pt]
        ({1.52*cos(205)},{1.52*sin(205)}) -- (-1.82,-.50);
      \node[anchor=east,inner sep=2pt] at (-1.84,-.50) {$\Gamma_D$};
      \draw[draw=sfgray,line width=.45pt]
        ({1.52*cos(17)},{1.52*sin(17)}) -- (1.75,.65);
      \node[anchor=west,inner sep=2pt] at (1.77,.65) {$\Gamma_N$};
      \foreach \ang in {-8,-28,-48}
        \draw[load] ({1.52*cos(\ang)},{1.52*sin(\ang)})
          -- ({1.96*cos(\ang)},{1.96*sin(\ang)});
      \node[anchor=west,inner sep=2pt] at (1.90,-.91)
        {$\bar{\mathbf t}$};
      \node[fill=sflight,inner sep=1.2pt] at (-.68,1.00) {$\Omega$};
    \end{scope}
  }

  \begin{scope}[shift={(3.00,6.80)}]
    \path[fill=white,draw=sfaccent,line width=.7pt]
      (-1.00,0) .. controls (-.48,.13) and (.48,.13) .. (1.00,0)
      .. controls (.48,-.13) and (-.48,-.13) .. cycle;
    \node[text=sfaccent,anchor=north,inner sep=1.5pt]
      at (-.92,-.14) {$J_{\mathbf u}$};
    \draw[load] (-.30,.19) -- (-.30,.72)
      node[right,midway,inner sep=2pt] {$\mathbf u^+$};
    \draw[load] (-.30,-.19) -- (-.30,-.72)
      node[right,midway,inner sep=2pt] {$\mathbf u^-$};
    \draw[load] (.64,.08) -- (.64,.69)
      node[right,near end,inner sep=2pt] {$\bm\nu_{\mathbf u}$};
  \end{scope}
  \node at (3.00,4.95)
    {$[\![\mathbf u]\!]=a\bm\nu_{\mathbf u},\quad a>0$};
  \node at (3.00,4.55) {(a)};

  \begin{scope}[shift={(10.30,6.80)}]
    \foreach \w/\extent/\tone in {.23/1.06/12,.14/1.01/30,.06/.95/70} {
      \path[fill=sfblue!\tone]
        (-\extent,0)
        .. controls (-.90,\w) and (-.55,\w) .. (0,\w)
        .. controls (.55,\w) and (.90,\w) .. (\extent,0)
        .. controls (.90,-\w) and (.55,-\w) .. (0,-\w)
        .. controls (-.55,-\w) and (-.90,-\w) .. cycle;
    }
    \node[anchor=east,fill=sflight,inner sep=1pt]
      at (-.48,.51) {$W^+$};
    \draw[draw=sfgray,line width=.45pt] (-.69,.38) -- (-.76,.15);
    \draw[black!80,dash pattern=on 2.4pt off 1.6pt,line width=.55pt]
      (.52,.18) circle[radius=.56];
    \fill[black!85] (.52,.18) circle[radius=.9pt];
    \node[anchor=south east,inner sep=2.2pt] at (.52,.18) {$\mathbf x$};
    \draw[load] (.52,.18) -- (.916,.576)
      node[midway,above left,inner sep=1.4pt] {$\varepsilon$};
    \node[fill=sflight,inner sep=1pt] at (.34,1.02)
      {$B_\varepsilon(\mathbf x)$};
    \node[note,fill=sflight,inner sep=1.4pt] at (0,-.95)
      {$W^-$: local};
  \end{scope}
  \node at (10.30,4.95)
    {$s_\varepsilon^+=\varepsilon\mathcal A_\varepsilon[W^+]$};
  \node at (10.30,4.55) {(b)};

  \begin{scope}[shift={(3.00,1.00)}]
    \foreach \xx in {1.8,3.6,5.4,7.2}
      \draw[black!16,densely dotted,line width=.3pt]
        (\xx,0) -- (\xx,2.50);
    \draw[black!16,densely dotted,line width=.3pt]
      (0,1.10) -- (7.60,1.10);
    \draw[draw=sfgray,dash pattern=on 3pt off 2pt,line width=.45pt]
      (0,2.20) -- (7.20,2.20);
    \node[anchor=east,inner sep=5pt] at (0,2.20) {$G_c/2$};
    \node[anchor=north east,note,inner sep=2pt,fill=white] at (7.20,1.88)
      {$2\Phi(\infty)=G_c$};
    \draw[draw=sfgray,dash pattern=on 3pt off 2pt,line width=.45pt]
      (0,0) -- (1.42,1.7356);
    \node[note,anchor=west,fill=white,inner sep=2pt]
      at (1.33,1.76) {$\Phi(s)\sim s$};
    \draw[sfaccent,line width=.75pt,domain=0:2,samples=180,smooth]
      plot ({3.6*\x},{2.2*(1-exp(-2*\x))});
    \draw[black!85,line width=.45pt] (0,0) rectangle (7.60,2.50);
    \node[right] at (7.60,0) {$s$};
    \node[above] at (0,2.50) {$\Phi(s)$};
    \foreach \xx/\lab in {1.8/{G_c/2},3.6/{G_c},5.4/{3G_c/2},7.2/{2G_c}} {
      \draw[black!85,line width=.45pt] (\xx,0) -- (\xx,.07);
      \draw[black!85,line width=.45pt] (\xx,2.50) -- (\xx,2.43);
      \node[anchor=north,inner sep=4pt] at (\xx,0) {$\lab$};
    }
    \foreach \yy in {1.10,2.20} {
      \draw[black!85,line width=.45pt] (0,\yy) -- (.07,\yy);
      \draw[black!85,line width=.45pt] (7.60,\yy) -- (7.53,\yy);
    }
    \node[anchor=north east,inner sep=3pt] at (0,0) {$0$};
    \node[anchor=center,fill=white,inner sep=2pt] at (4.78,.91)
      {$\displaystyle\Phi(s)=\frac{G_c}{2}\!\left(1-e^{-2s/G_c}\right)$};
  \end{scope}
  \node at (6.62,.22) {(c)};
\end{tikzpicture}
\endgroup

%% file: figures/dynamic_branching_setup.tex
\begin{tikzpicture}[
  x=1cm,y=1cm,
  font=\small,text=black,
  boundary/.style={draw=black!78,line width=0.65pt},
  extension/.style={draw=black!48,line width=0.35pt},
  dimension/.style={draw=black!65,line width=0.4pt,
    {Latex[length=1.6mm,width=1.0mm]}-{Latex[length=1.6mm,width=1.0mm]}},
  loading/.style={draw=blue!48!black,line width=0.8pt,
    -{Latex[length=2.1mm,width=1.5mm]}}
]
  \path[fill=black!4] (0,0) rectangle (8,3.2);
  \draw[boundary] (0,0) rectangle (8,3.2);

  \draw[draw=black,line width=1.1pt,line cap=butt] (0,1.6) -- (4,1.6);
  \node[anchor=north] at (2,1.43) {Initial slit $\Gamma_{0,\mathrm B}$};

  \foreach \x in {0.5,1.5,2.5,3.5,4.5,5.5,6.5,7.5} {
    \draw[loading] (\x,3.2) -- (\x,3.94);
    \draw[loading] (\x,0) -- (\x,-0.74);
  }
  \node[anchor=south] at (4,4.10)
    {$\overline{\mathbf t}=+\sigma_0\mathbf e_y,\quad \sigma_0=1\,\mathrm{MPa}$};
  \node[anchor=north] at (4,-0.88)
    {$\overline{\mathbf t}=-\sigma_0\mathbf e_y$};

  \draw[extension] (8.06,1.6) -- (8.45,1.6);
  \node[anchor=west] at (8.53,1.6) {Traction-free};

  \foreach \y in {0,1.6,3.2}
    \draw[extension] (-0.10,\y) -- (-1.0,\y);
  \draw[dimension] (-0.78,0) -- (-0.78,1.6)
    node[midway,fill=white,inner sep=2pt,rotate=90] {$20$};
  \draw[dimension] (-0.78,1.6) -- (-0.78,3.2)
    node[midway,fill=white,inner sep=2pt,rotate=90] {$20$};

  \foreach \x in {0,8}
    \draw[extension] (\x,-0.10) -- (\x,-1.99);
  \draw[extension] (4,-1.40) -- (4,-1.99);
  \draw[dimension] (0,-1.78) -- (4,-1.78)
    node[midway,fill=white,inner sep=2pt] {$50$};
  \draw[dimension] (4,-1.78) -- (8,-1.78)
    node[midway,fill=white,inner sep=2pt] {$50$};

  \node[anchor=west,text=black!65,font=\footnotesize]
    at (0,-2.32) {Dimensions in mm};
  \node[anchor=east] at (8,-2.32) {$t_p=1\,\mathrm{mm}$};
\end{tikzpicture}

%% file: figures/kalthoff_winkler_setup.tex
\begin{tikzpicture}[
  x=1cm,y=1cm,
  font=\small,text=black,
  boundary/.style={draw=black!78,line width=0.65pt},
  extension/.style={draw=black!48,line width=0.35pt},
  dimension/.style={draw=black!65,line width=0.4pt,
    {Latex[length=1.6mm,width=1.0mm]}-{Latex[length=1.6mm,width=1.0mm]}},
  loading/.style={draw=blue!48!black,line width=0.8pt,
    -{Latex[length=2.1mm,width=1.5mm]}}
]
  \path[fill=black!4] (0,0) rectangle (8,4);
  \draw[boundary] (0,0) rectangle (8,4);

  \draw[draw=black,line width=1.1pt,line cap=butt] (3,4) -- (3,2);
  \draw[draw=black,line width=1.1pt,line cap=butt] (5,4) -- (5,2);
  \node[anchor=east] at (2.75,3) {$\Gamma_{0,\mathrm l}$};
  \node[anchor=west] at (5.25,3) {$\Gamma_{0,\mathrm r}$};

  \draw[draw=blue!48!black,line width=1.1pt,line cap=butt] (3,4) -- (5,4);
  \foreach \x in {3.20,3.60,4.00,4.40,4.80}
    \draw[loading] (\x,4.72) -- (\x,4);
  \node[anchor=south,align=center,font=\footnotesize] at (4,4.90)
    {$u_y=-U(t),\, u_x$ free, 
     Linear velocity ramp\\
     $v_0=16.5\,\mathrm{m/s},\quad t_r=1\,\mu\mathrm{s}$};

  \draw[extension] (8.06,2) -- (8.35,2);
  \node[anchor=west,align=left,font=\footnotesize] at (8.43,2)
    {All other boundaries\\traction-free};

  \foreach \x in {0,3,5,8}
    \draw[extension] (\x,-0.10) -- (\x,-0.76);
  \draw[dimension] (0,-0.55) -- (3,-0.55)
    node[midway,fill=white,inner sep=2pt] {$75$};
  \draw[dimension] (3,-0.55) -- (5,-0.55)
    node[midway,fill=white,inner sep=2pt] {$50$};
  \draw[dimension] (5,-0.55) -- (8,-0.55)
    node[midway,fill=white,inner sep=2pt] {$75$};

  \foreach \y in {0,2,4}
    \draw[extension] (-0.10,\y) -- (-0.80,\y);
  \draw[dimension] (-0.58,0) -- (-0.58,2)
    node[midway,fill=white,inner sep=2pt] {$50$};
  \draw[dimension] (-0.58,2) -- (-0.58,4)
    node[midway,fill=white,inner sep=2pt] {$50$};
  \node[anchor=west,text=black!65,font=\footnotesize]
    at (0,-1.06) {Dimensions in mm};
  \node[anchor=east,font=\footnotesize]
    at (8,-1.06) {$t_p=1\,\mathrm{mm}$};
\end{tikzpicture}